\documentclass[a4paper,11pt]{amsart}
\usepackage[a4paper,hmargin=27mm,vmargin=28mm]{geometry}
\usepackage{amsmath,amssymb,mathtools,booktabs,array,enumitem}
\usepackage{atableau}
\usepackage{xcolor}
\usepackage{aliascnt}
\usepackage[
    colorlinks=true,
    linkcolor=blue!55!black,
    citecolor=blue!55!black,
    urlcolor=blue!70!black
]{hyperref}
\usepackage[nameinlink,capitalise,noabbrev]{cleveref}
\hypersetup{
    pdftitle={Scopes arrangements and product formulas},
    pdfauthor={Huang Lin and Tao Qin},
    pdfsubject={Scopes class products, face-centred bands and normalized abaci}
}
\setlist[enumerate]{label=\upshape(\alph*),leftmargin=*,itemsep=2pt}
\numberwithin{equation}{section}

\newaliascnt{theorem}{equation}
\newtheorem{theorem}[theorem]{Theorem}
\aliascntresetthe{theorem}

\newaliascnt{proposition}{equation}
\newtheorem{proposition}[proposition]{Proposition}
\aliascntresetthe{proposition}

\newaliascnt{lemma}{equation}
\newtheorem{lemma}[lemma]{Lemma}
\aliascntresetthe{lemma}

\newaliascnt{corollary}{equation}
\newtheorem{corollary}[corollary]{Corollary}
\aliascntresetthe{corollary}

\theoremstyle{definition}

\newaliascnt{definition}{equation}
\newtheorem{definition}[definition]{Definition}
\aliascntresetthe{definition}

\newaliascnt{example}{equation}
\newtheorem{example}[example]{Example}
\aliascntresetthe{example}

\newaliascnt{remark}{equation}
\newtheorem{remark}[remark]{Remark}
\aliascntresetthe{remark}

\newcommand{\Z}{\mathbb Z}

\newcommand{\R}{\mathbb R}
\newcommand{\Pcoweight}{\mathsf P^\vee}
\newcommand{\Qcoroot}{\mathsf Q^\vee}
\newcommand{\hcox}{h_{\mathring{\Phi}}}
\newcommand{\eps}{\varepsilon}
\newcommand{\blam}{\boldsymbol\lambda}
\newcommand{\bmu}{\boldsymbol\mu}
\newcommand{\cont}{\operatorname{cont}}
\newcommand{\defect}{\operatorname{def}}

\newcommand{\norm}[1]{\left\|#1\right\|}
\title[Scopes arrangements and product formulas]{Scopes arrangements and product formulas}
\author{Huang Lin}
\author{Tao Qin}
\date{}
\begin{document}
\begin{abstract}
    We develop normalized abacus combinatorics for affine type \(C\). We also give explicit product formulas for the number of Scopes classes at every level and every positive depth in affine types \(A\) and \(C\).
\end{abstract}
\maketitle

\section{Introduction}

The abacus is a useful combinatorial tool for studying Scopes equivalences of cyclotomic Hecke algebras in affine type \(A\). Affine categorical actions provide a general framework for these equivalences, and Webster describes the corresponding Scopes classes in terms of chambers of affine hyperplane arrangements \cite{Scopes,CR,Web23}. In this paper, we develop normalized abacus combinatorics for affine type \(C\) with arbitrary charge and apply it to the study of weight support, Scopes classes, and Rouquier weights.

Our abacus model pairs the \(2e\) runners of an ordinary abacus.
The resulting coordinates give explicit formulas for folded contents,
delta-depth and weight support at arbitrary charge and level
(\cref{thm:closed-minimum}). The support formula expresses classical
maximal-weight theory in these coordinates
\cite[Chapter~12]{K90}; related maximal-weight results in type \(C\)
appear in \cite{AP,CH,AHSW}.
It yields an exact criterion for Scopes equivalence
(\cref{thm:scopes}). At level one, the resulting reduction has a
unique normal form in each Scopes class
(\cref{thm:normal-forms}).

Our main counting result gives explicit product formulas at every
level and every positive depth in affine types \(A\) and \(C\).
Let \(\Lambda\) be a dominant integral weight of positive level,
and let \(\nu\in P(\Lambda)\) be dominant and delta-maximal,
meaning that \(\nu+\delta\notin P(\Lambda)\).
Write \(r\) and \(\hcox\) for the rank and Coxeter number of the
underlying finite root system, \(W\) for the affine Weyl group,
and \(W_\nu\) for the finite stabilizer of \(\nu\).
For every \(M\ge1\), let \(\chi_M\) be the characteristic polynomial
of the Scopes arrangement at \(\nu-M\delta\).
We construct an explicit multiset \(\mathcal E_{\Lambda,\nu}\)
of \(r\) integers, independent of \(M\), such that
\[
\begin{aligned}
 \chi_M(q)
   &=\prod_{b\in\mathcal E_{\Lambda,\nu}}(q-M\hcox-b),\\
 \#\mathrm{Scopes}(W(\nu-M\delta))
   &=\frac{1}{|W_\nu|}
       \prod_{b\in\mathcal E_{\Lambda,\nu}}(M\hcox+b+1).
\end{aligned}
\]
For each dominant maximal weight $\nu$, the multiset is determined by the corresponding content and the
coordinate blocks of \(\nu\)
(\cref{thm:all-level-product,thm:typeA-all-level-factorization}).
In particular, increasing depth translates every root of
\(\chi_M\) by the Coxeter number.
Chuang--Tan classify Scopes classes for Ariki--Koike algebras and
give a product counting RoCK classes \cite{CT}; our type-\(A\)
formula evaluates the total number of Scopes classes.

We also describe Scopes walls in every untwisted affine type.
They form a band centred at the normalized dominant weight,
together with two finite families translated with depth
(\cref{thm:band-decomposition}).
For highest-weight orbits, only the band remains, and its
characteristic polynomial is determined by the stabilizer of the
alcove face (\cref{thm:affine-band-factorization}).
The proof uses classical lattice-point and alcove-counting methods
\cite{Haiman,AthanasiadisCatalan,At05,Sommers}.
For general maximal weights, integral factorization already fails
in type \(G_2^{(1)}\) (\cref{prop:G2-counterexample}).

Finally, we prove that every abacus Rouquier weight in type \(C\)
is RoCK and that every RoCK Scopes class contains such a weight
(\cref{thm:rouquier-rock}). This is a type-\(C\) counterpart of
\cite[Proposition~4.6]{Web23}.
We also give an effective test for the Rouquier condition in terms
of content coordinates (\cref{thm:effective-rouquier}).

\section{Partitions and normalized abaci}\label{sec:abaci}

We first fix the normalized abacus coordinates and prove the level-one
content formulas. These coordinates also give the arbitrary-level
support calculation in \cref{sec:higher}.

\subsection{Residues and root conventions}
Fix \(e\ge2\) and put \(N=2e\). A partition is a weakly decreasing
sequence \(\lambda=(\lambda_1,\lambda_2,\ldots)\) of nonnegative integers
with finite sum \(|\lambda|\). Its length is the number \(\ell(\lambda)\)
of nonzero parts, and its Young diagram is
\[
 [\lambda]=\{(a,b):a\ge1,\ 1\le b\le\lambda_a\}.
\]
Rows and columns are numbered from the top and the left, respectively.
For \(x\in\Z\), let \(\bar x\in\{0,\ldots,N-1\}\) be its residue
modulo \(N\), and define \(\iota(x)=\min(\bar x,N-\bar x)\).
At integer charge \(t\), a node \((a,b)\) has integer content \(t+b-a\),
unfolded residue \(\overline{t+b-a}\), and folded residue
\(\iota(t+b-a)\). Set
\[
 c_i(\lambda;t)=\#\{(a,b)\in[\lambda]:\iota(t+b-a)=i\},
 \qquad \cont_t(\lambda)=\sum_{i=0}^e c_i(\lambda;t)\alpha_i.
\]
We suppress \(t\) and \(\lambda\) when the charge and partition are fixed.

Let \(\mathring{\mathfrak{h}}^*_{\mathbb{R}}=\R^e\) have orthonormal basis \(\eps_1,\ldots,\eps_e\).
Adjoin vectors \(\delta,\Lambda_0\) to get $\mathfrak{h}^*_{\mathbb{R}}$ and with extended non-degenerate symmetric bilinear form:
\[
 (\delta,\delta)=(\Lambda_0,\Lambda_0)=0,\quad
 (\delta,\mathring{\mathfrak{h}}^*_{\mathbb{R}})=(\Lambda_0,\mathring{\mathfrak{h}}^*_{\mathbb{R}})=0,\quad (\Lambda_0,\delta)=2.
\]
Then the coordinates of simple roots in type \(C_e^{(1)}\) Cartan datum are the following:
\begin{equation}\label{eq:roots}
 \alpha_0=\delta-2\eps_1,\quad
 \alpha_i=\eps_i-\eps_{i+1}\ (0<i<e),\quad
 \alpha_e=2\eps_e.
\end{equation}
Thus
\[
 \delta=\alpha_0+2\alpha_1+\cdots+2\alpha_{e-1}+\alpha_e,
 \qquad \operatorname{ht}(\delta)=N.
\]
Here \(\operatorname{ht}(\sum_i z_i\alpha_i)=\sum_i z_i\).
Write \(Q_+=\sum_i\Z_{\ge0}\alpha_i\). For \(0\le k\le e\), put
\[
 u(k)=\eps_1+\cdots+\eps_k,\qquad \Lambda_k=\Lambda_0+u(k).
\]
The empty sum is zero. The coroot pairing is
\(\langle \alpha^\vee_i,\nu\rangle=2(\alpha_i,\nu)/(\alpha_i,\alpha_i)\),
so \(\langle \alpha^\vee_i,\Lambda_k\rangle=\delta_{ik}\).
Each \(\Lambda_k\) has level \((\Lambda_k,\delta)/2=1\).
The symmetrizing integers are \(d_0=d_e=2\) and \(d_i=1\) for
\(0<i<e\). Our form is twice Kac's normalized invariant form
\cite[Chapters~6--7]{K90}; long roots therefore have squared length four.
The simple coroot pairings are unchanged by this rescaling.
For a dominant weight \(\Lambda\), define
\begin{equation}\label{eq:defect}
 \defect_\Lambda(\beta)=(\Lambda,\beta)-\tfrac12(\beta,\beta).
\end{equation}

\subsection{Beta numbers and paired runners}
The infinite charge-\(t\) beta set is
\[
 B_t(\lambda)=\{\lambda_a+t-a:a\ge1\}.
\]
It can be represented by a strictly decreasing sequence $(\beta_1,\beta_2,\ldots)$ so that $\beta_k = t-k$ for $k\gg 0$.
For $r\in [0,N)$, let $B^{(r)}_t(\lambda) := \{ x\in \mathbb{Z} : Nx+r \in B_t(\lambda) \}$. 
Obviously, $B^{(r)}_t(\lambda)$ is again a $\beta$-number of a partition $\lambda^{(r)}$ with charge $b_r$.
The $N$-partition $(\lambda^{(0)},\ldots ,\lambda^{(N-1)})$ is called the $N$-quotient of $\lambda$.
We represent a subset $B\subset \mathbb{Z}$ by placing beads at the positions indexed by the elements of $B$ on a runner. Then $(B_t^{(0)}(\lambda),\ldots ,B^{(N-1)}_t(\lambda))$ is called the abacus displayment of $\lambda$.

We pair physical runners by \(r\leftrightarrow N-1-r\), and define
\begin{equation}\label{eq:paired}
 p_i=b_{i-1}-b_{N-i},\qquad s_i=b_{i-1}+b_{N-i}\quad(1\le i\le e).
\end{equation}
Write \(p=(p_i)\), \(s=(s_i)\), and
\(h(p)=\#\{i:p_i\text{ is odd}\}\). Then
\begin{equation}\label{eq:parity}
 \sum_i s_i=t,\qquad s_i\equiv p_i\pmod2,\qquad
 b_{i-1}=\frac{s_i+p_i}{2},\quad b_{N-i}=\frac{s_i-p_i}{2}.
\end{equation}
The physical pairing differs from the residue folding
\(j\leftrightarrow-j\): adding a node moves a bead across the edge from
runner \(j-1\) to runner \(j\).

A runner $B^{(r)}$ is \emph{flush} if $b\in B^{(r)} \Rightarrow b-1\in B^{(r)}$, equivalently, $B^{(r)} = \{x\in \mathbb{Z} : x< b_r \}$.

\begin{definition}\label{def:weight}
An elementary upward slide moves a bead from \(x\) to an empty position
\(x-N\) on the same runner. The \emph{delta-weight}
\(w_\delta(\lambda)\) is the number of such slides required to flush
the abacus. Denote the resulting ordinary \(N\)-core by
\(\lambda^\circ\).
\end{definition}

This number and the terminal partition do not depend on the order of
slides. If the quotient positions on runner \(r\) are
\(q_{r,1}<\cdots<q_{r,n_r}\), then
\begin{equation}\label{eq:weight}
 w_\delta(\lambda)=\sum_{r=0}^{N-1}\sum_{a=1}^{n_r}(q_{r,a}-a+1).
\end{equation}
Each slide decreases this nonnegative integer by one. The beta-set
rim-hook correspondence \cite[Section~2.7]{JK} shows that a slide removes
a rim \(N\)-hook. Its nodes have consecutive integer contents and hence
folded content \(\delta\). Therefore
\begin{equation}\label{eq:flush}
 \cont_t(\lambda)=\cont_t(\lambda^\circ)+w_\delta(\lambda)\delta,
 \qquad |\lambda|=|\lambda^\circ|+Nw_\delta(\lambda).
\end{equation}

\subsection{Content coordinates}
We first work at a normalized charge \(0\le k\le e\).

\begin{proposition}\label{prop:coordinates}
Let \(\beta=\cont_k(\lambda)=\sum_i c_i\alpha_i\). Then
\begin{align}
 & c_0=w_\delta(\lambda)+\sum_r\binom{b_r}{2}
 =w_\delta(\lambda)+\frac{\norm{s}^2+\norm{p}^2-2k}{4},\label{eq:c0}\\
 & \beta=c_0\delta+u(k)-p,\label{eq:content}\\
 & c_i=2c_0+\min(i,k)-\sum_{j=1}^i p_j\quad(0<i<e),\qquad
 c_e=c_0+\frac{k-\sum_jp_j}{2},\label{eq:counts}\\
 & \defect_{\Lambda_k}(\beta)=2c_0+\frac{k-\norm{p}^2}{2}.
 \label{eq:defect-p}
\end{align}
In \eqref{eq:content}, vectors in \(\Z^e\) are identified with their
coordinates in the basis \((\eps_i)\), and the equality takes place in
\(E\oplus\R\delta\).
\end{proposition}
\begin{proof}
The contents in row \(a\) are \(k+1-a,\ldots,k+\lambda_a-a\).
The $0$-nodes are exactly the multiples of \(N\) in this
interval. Their number is
\[
 \#\{m\in\Z:k-a<Nm\le k+\lambda_a-a\}
 =\left\lfloor\frac{k+\lambda_a-a}{N}\right\rfloor
  -\left\lfloor\frac{k-a}{N}\right\rfloor.
\]
For a flush abacus, take $J\in \mathbb{Z}$ such that $L=NJ+k\geq  \ell (\lambda)$, we have
\begin{align*}
 c_0= & \sum_{r=1}^{\ell (\lambda)} \left\lfloor\frac{k+\lambda_a-a}{N}\right\rfloor
  -\left\lfloor\frac{k-a}{N}\right\rfloor \\
  =&\sum_{a=1}^{L}
 \left(\left\lfloor\frac{\lambda_a+L-a}{N}\right\rfloor   
 -\left\lfloor\frac{L-a}{N}\right\rfloor\right)\\
 = & \sum_{r=0}^{N-1}\left(\sum_{t=0}^{b_r+J}t - \sum_{t=0}^{J-\delta_{r\geq  k}}t\right) \\
 =& \sum_r\binom{J+b_r}{2}-k\binom{J+1}{2}-(N-k)\binom J2 \\
 =& \sum_{r< k}(b_r-1)J + \sum_{r\geq  k} b_rJ + \sum_r\binom{b_r}{2}  \\
 = & \sum_r\binom{b_r}{2}.
\end{align*}
where the last equality follows from \(\sum_r b_r=k\). Equation~\eqref{eq:flush} proves the first equality in \eqref{eq:c0} for arbitrary $\lambda$.
The second follows from \eqref{eq:parity}.

Adding a node replaces a beta number \(x-1\) by \(x\). If
\(j=\bar x\), the changes in the paired vector are
\[
\begin{array}{c|ccc}
j&0&i\text{ or }N-i\ (0<i<e)&e\\\hline
\Delta p&2\eps_1&-\eps_i+\eps_{i+1}&-2\eps_e.
\end{array}
\]
These are the negatives of the finite parts of the simple roots in
\eqref{eq:roots}. 
The empty partition has paired vector \(u(k)\).
Then \eqref{eq:content} and \eqref{eq:counts} follow from an inductive argument on $|\lambda|$. Finally,
\[
 (\Lambda_k,\beta)=2c_0+(u(k),u(k)-p),\quad
 (\beta,\beta)=\norm{u(k)-p}^2,\quad \norm{u(k)}^2=k.
\]
Substitution in \eqref{eq:defect} proves \eqref{eq:defect-p}.
\end{proof}

\begin{example}\label{ex:87}
Take \(e=2\), \(k=0\), and \(\lambda=(8,7)\). The following diagram
is filled with folded residues. We use \(J=1\), so \(L=4\), for the
abacus; its bead labels are beta numbers.
\begin{center}
\begin{tabular}{cc}
\Tableau{01210121,1012101}
&\Abacus[runner labels={0,1,2,3},entries=beads]{4}{8,7,0,0}
\end{tabular}
\end{center}
The beta numbers are \(11,9,1,0\). Thus
\[
 b=(0,1,-1,0),\quad p=(0,2),\quad s=(0,0),\quad
 w_\delta=3,\quad (c_0,c_1,c_2)=(4,8,3).
\]
In the second row, the integer contents are \(-1,0,1,2,3,4,5\).
The number of $0$-nodes is
\[
 \left\lfloor\frac{7+4-2}{4}\right\rfloor
 -\left\lfloor\frac{4-2}{4}\right\rfloor=2.
\]
Flushing gives \(\lambda^\circ=(2,1)\), of content
\(\alpha_0+2\alpha_1\).
\end{example}

\begin{corollary}\label{cor:content-p}
At fixed charge, two partitions have the same paired vector if and only
if their contents differ by an integer multiple of \(\delta\). They
have equal content if and only if they have equal paired vectors and
equal sizes.
\end{corollary}
\begin{proof}
Use \eqref{eq:content} and take heights; \(\operatorname{ht}(\delta)=N\).
\end{proof}

\begin{lemma}\label{lem:charges}
The operations
\[
 (\lambda,t)\longmapsto(\lambda,t+Na),\qquad
 (\lambda,t)\longmapsto(\lambda^{\mathrm t},-t)
 \quad(a\in\Z)
\]
preserve folded content, \(p\), and delta-weight. They reduce an
arbitrary charge to \(k=\iota(t)\). Equations~\eqref{eq:content} and
\eqref{eq:defect-p} hold for arbitrary \(t\) with \(k=\iota(t)\).
\end{lemma}
\begin{proof}
Translation by \(Na\) replaces every runner charge by \(b_r+a\).
Transposition gives the beta-set identity
\[
 B_{-t}(\lambda^{\mathrm t})
 =\Z\setminus\{-1-x:x\in B_t(\lambda)\},
\]
as follows by reversing the boundary of the diagram and exchanging its
horizontal and vertical steps. Hence transposition sends
\(b_r\) to \(-b_{N-1-r}\). Both operations fix \(p\); their effects
on the paired sums are \(s\mapsto s+2a\sum_{i=1}^e\varepsilon_i\) and \(s\mapsto-s\) respectively.
For a transposed node, the new integer content is the negative of
its old content; \(\iota(x)=\iota(-x)\). Charge translation adds
\(Na\) to every integer content. Thus both operations preserve the
entire folded content. In the beta-set identity, the map
\(Nq+r\mapsto-1-(Nq+r)=N(-q-1)+(N-1-r)\) exchanges beads and gaps,
which explains the cutoff change \(b_r(\lambda^{\mathrm t} ,-t)=-b_{N-1-r}(\lambda ,t)\).
Transposition preserves ordinary rim-hook weight, and charge translation
only translates its display.
Reduce \(t\) modulo \(N\), and if the result exceeds \(e\), negate
it and transpose. The remaining claims follow from
\cref{prop:coordinates}.
\end{proof}
For arbitrary \(t\), the first formula in \eqref{eq:c0} requires the
vacuum correction
\[
 c_0(\lambda;t)=w_\delta(\lambda)+
 \sum_r\left(\binom{b_r}{2}-\binom{v_r(t)}{2}\right),
 \qquad v_r(t)=\left\lceil\frac{t-r}{N}\right\rceil.
\]
Indeed \(v_r(t)\) is the cutoff of runner \(r\) for \(B_t(\varnothing)\),
and the same finite-display calculation proves the formula.

\section{Delta-maximal contents and their fibres}\label{sec:cores}

Fix \(0\le k\le e\) throughout this section.

\begin{definition}\label{def:depth}
The \emph{delta-depth} of a partition \(\lambda\) is
\[
 D_k(\lambda)=\max\{d\in\Z_{\ge0}:\cont_k(\lambda)-d\delta
                      =\cont_k(\mu)\text{ for some partition }\mu\}.
\]
A partition is a \emph{delta-core} if its delta-depth is zero.
\end{definition}
The maximum exists, since zero is allowed and \(d\le c_0(\lambda;k)\).
This definition does not require \([\mu]\subseteq[\lambda]\); see \cref{ex:52}.  

\subsection{The minimum at a fixed paired vector}
Put
\[
 \mathcal P_k=\{p\in\Z^e:\sum_i p_i\equiv k\pmod2\},\qquad
 \mathcal S_k(p)=\{s\in\Z^e:\sum_i s_i=k,\ s_i\equiv p_i\pmod2\}.
\]
The set \(\mathcal S_k(p)\) is nonempty exactly when
\(p\in\mathcal P_k\).

\begin{lemma}\label{lem:minimum}
Let \(p\in\mathcal P_k\) and \(h=h(p)\). Then
\begin{equation}\label{eq:min-square}
 \min_{s\in\mathcal S_k(p)}\norm{s}^2=k+|h-k|.
\end{equation}
The minimizing vectors have the following form:
\begin{enumerate}
\item If \(h\ge k\), the even positions of \(s\) have entries equal to zero;
the odd positions have entries equal to \(\pm 1\) with \((h-k)/2\) negative entries.
\item If \(h\le k\), the odd positions have entries equal to \(1\); the even
positions have entries \(0,2\), with \((k-h)/2\) entries equal to \(2\).
\end{enumerate}
Equivalently, \(s\) is minimizing if and only if its range
\(\max_i s_i-\min_i s_i\le2\).
\end{lemma}
\begin{proof}
It is obvious that $h\equiv k\pmod 2$.  For \(h\ge k\), each odd entry contributes at least one to the
squared norm, giving the bound \(h\). Equality forces precisely the
entries in (a), and the sum fixes the number of negative entries.
For \(h\le k\), use
\[
 x^2\ge2x-1\quad(x\text{ odd}),\qquad
 x^2\ge2x\quad(x\text{ even}).
\]
The equality holds only at \(1\), respectively at \(0,2\). Hence 
\(\norm{s}^2 \geq  2k-h\), and (b) follows.

The displayed minimizers have range at most \(2\). 
Conversely, for a feasible vector of range at most two, \(\min_is_i \leq  -2\) (resp. \(\max_is_i\geq  3\)) would force \(\sum_is_i<0\) (resp.\(\sum_is_i> e\)). 
Therefore, its entries lie in \(\{-1,0,1,2\}\), and \(-1\) and \(2\) cannot both occur. 
The parity and sum condition then give (a) or (b).
\end{proof}

\begin{theorem}\label{thm:depth}
For \(p\in\mathcal P_k\), set the cost of $p$ to be
\[
 f_k(p)=\frac{\norm{p}^2-k+|h(p)-k|}{4}.
\]
Among partitions with paired vector \(p\), we have \(c_0 \geqslant f_k(p)\). Moreover,
\begin{align}
 &D_k(\lambda)=c_0-f_k(p)
 =w_\delta(\lambda)+\frac{\norm{s}^2-k-|h(p)-k|}{4},\label{eq:depth}\\
 &\defect_{\Lambda_k}(\cont_k(\lambda))=2D_k(\lambda)+\frac{|h(p)-k|}{2}.
 \label{eq:defect-depth}
\end{align}
A partition is a delta-core if and only if its abacus is flush and its
paired sums have range at most two.
\end{theorem}
\begin{proof}
Equation~\eqref{eq:c0} and \cref{lem:minimum} give the bound
\(c_0\ge f_k(p)\) and the second equality in \eqref{eq:depth}. Choose a minimizing \(s\), reconstruct the
integers \(b_r\) using \eqref{eq:parity}, and take the flush abacus
with these cutoffs. Its total charge is \(k\), so it decodes to a
partition $\mu$. It attains \(c_0(\mu) =f_k(p)\), also proving that this number is a nonnegative integer. 
By \eqref{eq:content}, $\cont_k(\lambda) -\cont_k(\lambda') \in\Z \delta$ if and only if $p(\lambda) =p (\lambda')$. 
We claim that $D_k(\mu) =0$ since otherwise $\cont_k(\mu) = \cont_k(\mu') + d\delta$ for some $d>0$ and $c_0(\mu')=c_0(\mu) -d <f_k(p)$ contradicting the above argument. 
Now $D_k(\lambda) =c_0(\lambda) -c_0(\mu) = c_0 -f_k(p)$.
This completes the proof of depth formula.

Substitution in
\eqref{eq:defect-p} gives \eqref{eq:defect-depth}. 
The last statement is a direct consequence of \eqref{eq:depth}.
\end{proof}

\subsection{Reduction and uniqueness}
There are two reduction operations.
\begin{enumerate}
\item \emph{Flushing.} Perform elementary upward slides until every runner is flush. 
The paired coordinates remain unchanged, and the difference of contents is \(w_\delta(\lambda)\) copies of \(\delta\).
\item \emph{Balancing.} Given an $N$-core partition $\lambda$, choose \(i,j\) with balancing gap
\(s_i-s_j>2\), and reconstruct an $N$-core partition $\mu$ satisfying
\begin{equation}\label{eq:balance}
 s(\mu) =s(\lambda)-2\varepsilon_i+2\varepsilon_j,\qquad p(\mu) =p(\lambda).
\end{equation}
\end{enumerate}

\begin{theorem}\label{thm:reduction}
Every sequence of reductions terminates at a delta-core \(\rho\), and
\[
 \cont_k(\rho)=\cont_k(\lambda)-D_k(\lambda)\delta.
\]
If the successive balancing gaps are \(t_1,\ldots,t_a\), then
\begin{equation*}
 D_k(\lambda)=w_\delta(\lambda)+\sum_{r=1}^a(t_r-2).
\end{equation*}
\end{theorem}
\begin{proof}
Note that the reductions preserve the coordinate $p$. 
The range of $s$ decreases after finitely many steps since $\sum_{i<j}|s_i -s_j|$ decreases at each balancing step. 
Therefore, every sequence of reductions terminates at $\rho$ and $D_k(\rho) =0$ by \cref{thm:depth}. 
Therefore, we have
\begin{align*}
    \cont_k(\lambda) - \cont_k(\rho) & \overset{\eqref{eq:content}}{=} \left(c_0(\lambda)-c_0(\rho)\right)\delta \\
    & \overset{\eqref{eq:depth}}{=} \left(D_k(\lambda)-D_k(\rho)\right)\delta \\
    &\; = D_k(\lambda)\delta.
\end{align*}
By \eqref{eq:depth}, $D_k(\lambda) = w_\delta (\lambda) + D_k(\lambda^\circ)$. 
If $N$-core $\mu'$ is obtained from $N$-core $\mu$ by a single balancing, say, $s_i(\mu') = s_i(\mu)-2$ and $s_j(\mu') =s_j(\mu)+2$, then, by \eqref{eq:depth}, $D_k(\mu)-D_k(\mu') = \norm{s(\mu)} - \norm{s(\mu')} = s_i(\mu) - s_j(\mu')-2$. 
This proves the last statement.
\end{proof}

Let \(\mathcal C_k\) be the set of delta-core partitions. Define
\(\rho\sim_k\sigma\) for $\rho,\sigma\in \mathcal{C}_k$ if $\cont_k(\rho) = \cont_k(\sigma)$. 
For any partition \(\lambda\), write \(\mathfrak c_k(\lambda)\) for the class
of a terminal core. \cref{thm:reduction} makes this well defined.

\begin{corollary}\label{cor:core-class}
The map \([\rho]\mapsto p(\rho)\) is a bijection
\(\mathcal C_k/\!\sim_k\ \longrightarrow\mathcal P_k\).
For partitions \(\lambda,\mu\), the following are equivalent:
\begin{enumerate}
\item \(\cont_k(\lambda)=\cont_k(\mu)\);
\item \(\mathfrak c_k(\lambda)=\mathfrak c_k(\mu)\) and
\(|\lambda|=|\mu|\);
\item \(\mathfrak c_k(\lambda)=\mathfrak c_k(\mu)\) and
\(D_k(\lambda)=D_k(\mu)\).
\end{enumerate}
\end{corollary}
\begin{proof}
The map is well-defined by \eqref{eq:content}. 
The injectivity is a consequence of \eqref{eq:content} and \eqref{eq:depth}.
The reconstruction in \cref{thm:depth} proves surjectivity.
For a representative \(\rho\) of the core class,
\[
 \cont_k(\lambda)=\cont_k(\rho)+D_k(\lambda)\delta,
 \qquad |\lambda|=|\rho|+ND_k(\lambda).
\]
The equivalences follow.
\end{proof}

\subsection{Representatives of a core class}
Fix \(p\in\mathcal P_k\). 
By \cref{lem:minimum}, there exists a minimizer $s\in \mathcal{S}_k(p)$ that is lexicographically smallest among all minimizers which we denote by $s^*(p,k)$. 
Let \(\rho (k,p)\) be the $N$-core partition reconstructed from \(p,s^*(p,k)\).

By \cref{thm:depth}, delta-core partitions with a fixed paired vector $p$ are in bijection with the set of minimizers $s\in \mathcal{S}_k(p)$.
Note that we can also apply \eqref{eq:balance} to a partition with balancing gap $s_i-s_j =2$.
These \emph{zero-cost balancings} permute minimizers and hence act transitively on each $\sim_k$-class.
We define a graph $J(k,p)$ whose vertices are the delta-core partitions with paired vector $p$, with an edge connecting $\rho$ and $\sigma$ if one can be obtained from the other by a zero-cost balancing operation.
Then the following proposition is clear.

\begin{proposition}\label{prop:core-fibre}
The number of delta-core partitions with paired vector
\(p\) is
\begin{equation}\label{eq:g}
 g_k(p)=
 \begin{cases}
 \displaystyle\binom{h(p)}{(h(p)-k)/2},&h(p)\ge k,\\[4pt]
 \displaystyle\binom{e-h(p)}{(k-h(p))/2},&h(p)\le k.
 \end{cases}
\end{equation}
The graph $J(k,p)$ is respectively
\(J(h(p),(h(p)-k)/2)\) or \(J(e-h(p),(k-h(p))/2)\), where
\(J(n,a)\) has the \(a\)-element subsets of $[n]$ as vertices and replaces one element at each edge. 
In particular, this graph is connected.
\end{proposition}

\begin{example}\label{ex:44}
Let \(e=3\), \(k=0\), and \(\lambda=(4,4)\). Then
\[
 b=(0,0,1,1,-1,-1),\quad p=(1,1,0),\quad s=(-1,-1,2),
 \qquad w_\delta=0,\quad D_0=1.
\]
The following are flush displays with \(L=6\) beads.
\begin{center}
\begin{tabular}{ccc}
\Abacus[runner labels={0,1,2,3,4,5},entries=beads]{6}{4,4,0,0,0,0}
&\Abacus[runner labels={0,1,2,3,4,5},entries=beads]{6}{2,0,0,0,0,0}
&\Abacus[runner labels={0,1,2,3,4,5},entries=beads]{6}{1,1,0,0,0,0}\\
\((4,4)\)&\((2)\)&\((1,1)\)
\end{tabular}
\end{center}
Balancing against the second, respectively first, coordinate gives
\(s=(-1,1,0)\), respectively \((1,-1,0)\). Both subtract one copy
of \(\delta\). The two core diagrams have the same content:
\begin{center}
\Tableau{01}\qquad\Tableau{0,1}.
\end{center}
A zero-cost move exchanges them. Our sorting convention selects
\(\rho(0,(1,1,0))=(2)\).
\end{example}

\begin{example}\label{ex:52}
For \(e=2\), \(k=1\), and \(\lambda=(5,2)\), one has
\(p=(1,2)\), \(s=(-1,2)\), and \(w_\delta=0\).
Balancing gives \(s=(1,0)\) and the unique core \((1,1,1)\).
\begin{center}
\begin{tabular}{ccc}
\Tableau{12101,01}&\(\longrightarrow\)&\Tableau{1,0,1}\\[4pt]
\Abacus[runner labels={0,1,2,3},entries=beads]{4}{5,2,0,0,0}
& &\Abacus[runner labels={0,1,2,3},entries=beads]{4}{1,1,1,0,0}
\end{tabular}
\end{center}
The displays have \(L=5\). Their contents are
\(2\alpha_0+4\alpha_1+\alpha_2\) and \(\alpha_0+2\alpha_1\), so
the depth is one. The core has three rows and is not a subdiagram of
\((5,2)\). Thus a reduction by box deletions alone cannot always work.
\end{example}

\begin{example}\label{ex:weight-depth}
At \(e=2,k=0\), the partitions \((8)\) and \((4,2,2)\) both have
content \(2\delta\), paired vector zero, and depth two. The first has
delta-weight two. The second is flush, with \(s=(2,-2)\), and has
delta-weight zero. One balancing move takes it to the empty partition
and removes \(2\delta\). Thus delta-depth counts copies of \(\delta\)
removed, not balancing iterations; delta-weight is not a content invariant.
\end{example}

\subsection{A generating series for a content string}
For completeness, the paired coordinates also describe every partition on a fixed delta-string. 
The ordinary \(N\)-quotient consists of one
partition on each runner, whose total size is \(w_\delta\). It can be
read from the displacements in \eqref{eq:weight}; see
\cite[Section~2.7]{JK}.

\begin{proposition}\label{prop:series}
For \(p\in\mathcal P_k\), the following is an identity of formal power
series in \(z\):
\[
 \sum_{\lambda:\,p(\lambda)=p}z^{D_k(\lambda)}
 =\frac{\displaystyle\sum_{s\in\mathcal S_k(p)}
                z^{(\norm{s}^2-k-|h(p)-k|)/4}}
        {\displaystyle\prod_{m\ge1}(1-z^m)^N}.
\]
Its constant coefficient is \(g_k(p)\).
\end{proposition}
\begin{proof}
We have
\begin{align*}
    \sum_{\lambda : p(\lambda) =p} z^{D_k(\lambda)} & \overset{\eqref{eq:depth}}{=} \sum_{\lambda : p(\lambda) =p}z^{w_\delta (\lambda)+(\norm{s(\lambda)}^2-k-|h(p)-k|)/4} \\
    & = \sum_{s\in\mathcal{S}_k(p)}\left(z^{(\norm{s}^2-k-|h(p)-k|)/4}\sum_{\lambda:p(\lambda)=p,s(\lambda)=s}z^{w_\delta (\lambda)}\right)
\end{align*}
Since $w_\delta (\lambda)$ is equal to the size of the $N$-quotient of $\lambda$, so $\sum_{\lambda : p(\lambda)=p,s(\lambda)=s}z^{w_\delta (\lambda)}$ is independent of $p$ and $s$ and equal to the generating series of $N$-partitions.
Therefore,
\begin{align*}
    \sum_{\lambda : p(\lambda) =p} z^{D_k(\lambda)} & = \sum_{s\in \mathcal{S}_k(p)}z^{(\norm{s}^2-k-|h(p)-k|)/4}\left(\prod_{m\ge 1}\frac{1}{(1-z^m)}\right)^N
\end{align*}
This proves the desired identity.
The constant term follows from
\cref{prop:core-fibre}.
\end{proof}

The numerator is a normalized coset theta series. Indeed, for any
\(s^0\in\mathcal S_k(p)\),
\[
 \mathcal S_k(p)=s^0+2A_{e-1},\qquad
 A_{e-1}=\{z\in\Z^e:\textstyle\sum_i z_i=0\}.
\]
It depends only on \(k\) and \(h(p)\), since a permutation carries
any parity pattern of that size to any other. Thus the number of labels
at fixed depth is constant on a Weyl orbit. This is a series for the
full partition Fock space. An irreducible string function instead
counts weight multiplicities in \(V(\Lambda_k)\); the two already
have different constant terms in the defect-one examples of
\cref{prop:zero-simples}.

\section{Multipartitions and minimum content}\label{sec:higher}

Fix an interger $\ell>0$ and a multicharge \(\mathbf k=(k_0,\ldots,k_{\ell-1})\) with
\(0\le k_c\le e\). 
An \(\ell\)-partition is an ordered tuple of partitions \(\blam=(\lambda^{(0)},\ldots,\lambda^{(\ell-1)})\).
Its size and content are the sums of the component sizes and contents.
Write
\begin{align*}
 & \Lambda=\sum_{c=0}^{\ell -1}\Lambda_{k_c},\quad u(\mathbf{k})=\sum_{c=0}^{\ell -1}u(k_c),\quad
 P=P(\blam,\mathbf{k})=\sum_{c=0}^{\ell -1}p^{(c)},\\ 
 & p^{(c)}=p(B_{k_c}(\lambda^{(c)})),\quad \cont_{\mathbf k}(\blam)=\sum_{c=0}^{\ell -1}\cont_{k_c}(\lambda^{(c)}).
\end{align*}
Summing \eqref{eq:content} and using \(\Lambda=\ell\Lambda_0+u\) gives
\begin{equation}\label{eq:higher-coordinates}
 \cont_{\mathbf{k}}(\blam)=c_0\delta+u-P,\qquad
 d:=\defect_\Lambda(\beta)=2\ell c_0+\frac{\norm{u}^2-\norm{P}^2}{2}.
\end{equation}
In particular, \(P\) and \(c_0\) determine all content. Explicitly,
\[
 c_i=2c_0+\sum_{r=1}^i(u_r-P_r)\ (0<i<e),\qquad
 c_e=c_0+\tfrac12\sum_r(u_r-P_r).
\]

We say that that \(P\in \mathbb{Z}^e\) is $\mathbf{k}$-admissible if \(\sum_iP_i\equiv\sum_ck_c\pmod2\).
An ($\mathbf{k}$-)\emph{allocation} of \(P\) is \((p^{(c)})\in \prod_{c=0}^{\ell}\mathcal{P}_{k_c}\) with \(\sum_cp^{(c)}=P\) which is denoted by $(p^{(c)})\vdash P$.
Such allocations exist, for example, $p^{(c)} = u(k_c)$ for $c<\ell -1$ and $p^{(\ell -1)}=P-\sum_{c=0}^{\ell -2}u (k_c)$.
Equivalently, the admissible sequence is the image of the sum map $\prod_{c=0}^{\ell -1}\mathcal{P}_{k_c}\to \mathbb{Z}^e$ and allocations of $P$ are preimages of $P$.
Define
\begin{equation}\label{eq:T}
 T_{\mathbf k}(P)=\min_{(p^{(c)})\vdash P}  \sum_cf_{k_c}(p^{(c)}).
\end{equation}
This is initially an effective finite minimum; \cref{thm:closed-minimum}
will evaluate it without a search. Given any allocation $(p^{(c)})\vdash P$ of cost \(C_*=\sum_{c=0}^{\ell -1}f_{k_c}(p^{(c)})\), 
every allocation of cost at most \(C_*\) satisfies
\begin{equation}\label{eq:finite-bound}
 \norm{p^{(c)}}^2\le4C_*+k_c.
\end{equation}
It suffices to enumerate the finitely many integer vectors in these balls.

\begin{definition}
    The \emph{delta-depth} of a charged $\ell$-partition $\blam$ is
    \[
    D_{\mathbf k}(\blam)=\max\{d\in \mathbb{Z}_{\ge 0}: \cont_{\mathbf{k}}(\blam) - d\delta =\cont_{\mathbf{k}}(\bmu) \text{ for some charged $\ell$-partition $\bmu$}\}.
    \]
\end{definition}

\begin{theorem}\label{thm:higher}
Let $(\blam ,\mathbf{k})$ be a charged $\ell$-partition with total paired vector $P$, then the number of $0$-nodes $c_0\ge T_{\mathbf{k}}(P)$. 
Moreover,
\begin{equation}\label{eq:higher-depth}
 D_{\mathbf k}(\blam)=c_0-T_{\mathbf k}(P)
 =\sum_cD_{k_c}(\lambda^{(c)})+
   \sum_cf_{k_c}(p^{(c)})-T_{\mathbf k}(P).
\end{equation}
An $\ell$-partition has depth zero if and only if every component is a
delta-core and its paired vectors form a minimizing allocation.
\end{theorem}
\begin{proof}
The last equality of \eqref{eq:higher-depth} follows from \cref{thm:depth} and hence $c_0\ge T_{\mathbf{k}}(P)$ by definition.
The equality holds if and only if $(p^{(c)})$ is minimal allocation and $\lambda^{(c)}$ is a delta-core with respect to $k_c$ for $0\le c<\ell$.
We claim that such $\ell$-partition $\blam$ has depth zero since otherwise $\cont_{\mathbf{k}}(\blam) =\cont_{\mathbf{k}}(\bmu) + d\delta$ for some $d>0$ and $c_0(\bmu) - T_{\mathbf{k}}(P) <0$, a contradiction.
By \eqref{eq:higher-coordinates}, we can complete the proof by showing the existence of depth zero $\ell$-partitions with total paired vector $P$.
Take a minimizing allocation $(\widetilde{p}^{(c)})\vdash P$ and a minimizing $s^{(c)}\in\mathcal{S}_{k_c}(\widetilde{p}^{(c)})$, we can reconstruct a charged $\ell$-partition $\bmu$ whose $c$-th component is the depth zero partition with paired vectors $\widetilde{p}^{(c)},s^{(c)}$.
Clearly, $c_0(\bmu) - T_{\mathbf{k}}(P) =0$ and hence $D_{\mathbf{k}}(\bmu) =0$.
\end{proof}

Again, for any two depth-zero $\ell$-partition, we define $\blam\sim_{\mathbf{k}}\bmu$ if $\cont_{\mathbf{k}}(\blam) = \cont_{\mathbf{k}}(\bmu)$.
Each admissible \(P\) determines exactly one class, of content
\(T_{\mathbf k}(P)\delta+u-P\). 
In general, content is classified by the core class and size, equivalently,  by the core class and depth.
The proof is the same coefficient-and-height argument as in
\cref{cor:core-class}.

In the proof of \cref{thm:higher}, we give reconstruction algorithms of depth-zero multi-partitions. 
The number of delta-zero multi-partitions is
\[
 \sum_{\substack{(p^{(c)})\vdash P\\
                  \sum_cf_{k_c}(p^{(c)})=T_{\mathbf k}(P)}}
                  \prod_cg_{k_c}(p^{(c)}).
\]
These counts concern partition labels, not simple modules. 
The multipartition delta-weight is \(\sum_cw_\delta(\lambda^{(c)})\);
it counts componentwise flushing only.

\begin{example}\label{ex:allocation}
Let \(e=2\), \(\mathbf k=(0,1)\), and \(\blam=((1),(3))\).
Its residue-filled component diagrams are
\begin{center}
\(\bigl(\,\)\Tableau{0}\(\,,\ \)\Tableau{121}\(\,\bigr)\).
\end{center}
The paired vectors are \((2,0)\) and \((-1,0)\). Each component is a
defect-zero delta-core, but their total content is \(\delta\).
The empty multipartition has the same total vector \(P=(1,0)\), so
\(T_{(0,1)}(P)=0\) and the depth is $1$. Its delta-weight is $0$.
\end{example}

\begin{proposition}\label{prop:support}
The weight $\Lambda -\beta \in P(\Lambda)$ if and only if $\beta = \cont_{\mathbf{k}}(\blam)$.
Consequently \eqref{eq:higher-depth} is also depth in this weight support.
\end{proposition}
\begin{proof}
Recall the Fock space representation $\mathcal{F}(\Lambda)$, see, for example, 
\cite[Corollaries~6D.4 and~6D.6]{EM24}. Its basis vector at \(\blam\)
has weight \(\Lambda-\cont_{\mathbf k}(\blam)\), and its empty vector
generates \(V(\Lambda)\). This gives one inclusion. 
Conversely, take a standard tableau $\mathsf{T}$ of $\blam$, we have 
\[
F_{\mathsf{T}(|\blam|)}\cdots F_{\mathsf{T}(1)}(\emptyset) = g(q)\blam + \sum_{\bmu\neq \blam}g'(q)\bmu \in V(\Lambda) \subset \mathcal{F}(\Lambda). 
\]
for some nonzero $g(q)\in \mathbb{N}[q,q^{-1}]$, so $\Lambda -\cont_{\mathbf{k}}(\blam) \in P(\Lambda)$. This completes the first part.
The last statement follows immediately.
\end{proof}
Alternatively, we can use the minimal categorification of $V(\Lambda)$ to prove this proposition. By \cite{EM24}, $\Lambda -\beta\in P(\Lambda)$ if and only if the cyclotomic KLR algebra $R^\Lambda_\beta \neq 0$ which is equivalent to the existence of a Specht module $S(\blam)$ with $\cont_{\mathbf{k}}(\blam) =\beta$.

This identifies support only; it asserts no equality of weight multiplicities
with the number of all multipartition labels.

\begin{theorem}\label{thm:closed-minimum}
Let \(P\) be admissible. Fold each coordinate modulo \(2\ell\) into
\([0,\ell]\), and sort the resulting coordinates decreasingly. Denote
the resulting vector by \(a=(a_1,\ldots,a_e)\), and put
\begin{equation*} 
 \zeta_\Lambda(a)=\left\lceil\frac12
       \max_{0\le j\le e}\sum_{i=1}^j(a_i-u_i)\right\rceil.
\end{equation*}
Then
\begin{equation}\label{eq:closed-T}
 T_{\mathbf k}(P)=
       \frac{\norm{P}^2-\norm{a}^2}{4\ell}+\zeta_\Lambda(a).
\end{equation}
In particular, the numerical minimum is computable in \(O(e\log e)\) 
arithmetic operations once \(u\) is given.
The dominant maximal weights of \(V(\Lambda)\) are exactly
\[
 \ell\Lambda_0+a-\zeta_\Lambda(a)\delta,
 \qquad \ell\ge a_1\ge\cdots\ge a_e\ge0,\quad
 \sum_i a_i\equiv\sum_c k_c\pmod2.
\]
\end{theorem}
\begin{proof}
Recall the natural action of $W =W(C^{(1)}_e)$ on $\mathfrak{h}^*_{\mathbb{R}}$. Note that the level coset
$$\mathfrak{h}^*_{\ell} = \{ x\in \mathfrak{h}^*_{\mathbb{R}}/\mathbb{R}\delta \mid (x,\delta)/2 =\ell \}\simeq \ell \Lambda_0 + \mathring{\mathfrak{h}}^*_{\mathbb{R}}$$
is stable under the induced $W$-action. 
Recall that $W$ has an action on $\mathring{\mathfrak{h}}^*_{\mathbb{R}}$ given by signed permutations and translations in \(2\ell\Z^e\) so that the identification $\mathfrak{h}^*_\ell \simeq \mathring{\mathfrak{h}}^*_{\mathbb{R}}$ is $W$-equivariant. 
The fundamental domain is the closed alcove \(\ell\ge x_1\ge \cdots\ge x_e\ge 0\). 
Thus \(a\) is the unique point of the closed alcove in \(W\cdot P\). 
Write $a =w(P)$, then $w(\ell \Lambda_0 + P-c_0\delta) =\ell \Lambda_0 + a -c'_0\delta$.
Since \(\norm{\ell\Lambda_0+P-c_0\delta} =\norm{P}^2-4\ell c_0\), we have
\[
 c'_0=c_0-\frac{\norm{P}^2-\norm{a}^2}{4\ell}.
\]
Combining \cref{thm:higher} and \cref{prop:support}, we see that \eqref{eq:closed-T} will follow from $T_{\mathbf{k}}(a) =  \zeta_{\Lambda}(a)$.

For $\nu \in P^+$, it is known that \(\nu \in P(\Lambda)\) if and only if \(\Lambda-\nu\in Q_+\)
\cite[Proposition~11.3(a)]{K90}; see also \cite[Lemma~3.6]{Web23}. 
Applied to
\(\nu=\ell\Lambda_0+a-c'_0\delta\) we obtain 
\[
 2c'_0\ge\sum_{i=1}^j(a_i-u_i)\quad(0\le j\le e).
\]
Therefore, $\nu$ is a dominant maximal weight of $V(\Lambda)$ if and only if $T_{\mathbf{k}}(a) =c'_0 =\zeta_{\Lambda}(a)$.
Sorting and taking partial sums gives the operation bound. 
\end{proof}

\cref{thm:closed-minimum} evaluates the allocation minimum; 
it does not give an equally efficient algorithm for constructing a
minimizing allocation. The finite reconstruction above remains valid.
The description of dominant maximal weights is a coordinate form of the
classical support criterion, to be compared with the explicit
constructions in \cite[Proposition~3.6]{AHSW}. It is not a new
classification of maximal weights.

\begin{corollary}\label{cor:weyl-depth}
Depth and defect are $W$-invariant. Moreover,
\[
 T_{\mathbf k}(2\ell-P_1,P_2,\ldots,P_e)
       =T_{\mathbf k}(P)+\ell-P_1,
\]
and \(T_{\mathbf k}\) is invariant under signed permutations.
\end{corollary}
\begin{proof}
    Using \cref{prop:support}, the first statement is classical, see \cite{K90}. The rest statements are direct consequences of \eqref{eq:closed-T}.
\end{proof}

\begin{corollary}\label{cor:maximal-defect}
Let \(P,a\in \mathbb{Z}^e\) be as in \cref{thm:closed-minimum}, put
\[
 d_\Lambda^{(0)}(a)=\frac{\norm u^2-\norm a^2}{2}+2\ell\zeta_\Lambda(a).
\]
Then,
\begin{equation*} 
 \defect_\Lambda (c_0\delta +u-P)=d_\Lambda^{(0)}(a) +2\ell\bigl(c_0-T_{\mathbf k}(P)\bigr).
\end{equation*}
Consequently, \(d_\Lambda^{(0)}(a)\) is the exact minimal defect
on the delta-string at fixed \(P\), attained at its maximal weight. 
\end{corollary}
\begin{proof}
Substitute \eqref{eq:closed-T} into \eqref{eq:higher-coordinates}.
All nonnegative depths occur by \cref{thm:higher}.
\end{proof}


\begin{example}\label{ex:closed-T}
Take \(e=3\), \(\mathbf k=(0,2)\), and \(P=(5,-2,1)\).
Here \(\ell=2\), \(u=(1,1,0)\), and \(a=(2,1,1)\).
The partial sums of \(a-u\), including the empty sum, are
\(0,1,1,2\). Thus \(\zeta_\Lambda(a)=1\) and
\[
 T_{(0,2)}(P)=\frac{30-6}{8}+1=4.
\]
For example, the allocation \(((2,-1,1),(3,-1,0))\) has costs
\(f_0(2,-1,1)=2\) and \(f_2(3,-1,0)=2\), so it attains the
minimum without an allocation search.
\end{example}

\section{The Scopes criterion}\label{sec:scopes}

\subsection{The categorical setting}
Fix an algebraically closed field \(\Bbbk\), of arbitrary characteristic.
Let \(R^\Lambda(\beta)\) be the cyclotomic KLR algebra for the root
datum \eqref{eq:roots}.
The generators, relations and cyclotomic
quotient are those of \cite[Definition~2C.2 and (2C.3)]{EM24}; explicitly
the latter relation is
\(y_1^{\langle \alpha^\vee_{\nu_1},\Lambda\rangle}e(\nu)=0\).
Here \(e(\nu)\) is the idempotent for a word whose simple roots sum to
\(\beta\), and \(y_1\) is its first polynomial generator. The grading
is \(\deg y_re(\nu)=2d_{\nu_r}\) and
\(\deg\psi_re(\nu)=-(\alpha_{\nu_r},\alpha_{\nu_{r+1}})\).
For \(\beta\notin Q_+\), set \(R^\Lambda(\beta)=0\).

Write \(\mathcal C_\mu=R^\Lambda(\beta)\text{-}\mathrm{gmod}\), where
\(\mu=\Lambda-\beta\), for finite-dimensional modules, with grading
forgotten. 
These finite-dimensional algebras are the \emph{weight summands} of \(R^\Lambda(n)\). 
We use this term without an indecomposability assertion; compare the discussion after \cite[Definition~4A.1]{EM24}. 
A classification of these summands by content does not itself classify their indecomposable blocks. 

For $i\in \{ 0,1,\ldots ,e\}$, there is a pair of functors \( (F_i,E_i)\) called the $i$-induction and $i$-restriction functor respectively, they are denoted by $F_i^\Lambda,E_i^\Lambda$ in \cite{KK}.
Up to a graded shift, $F_i$ is biadjoint to $E_i$.
We have divided powers $F_i^{(a)}, E_i^{(a)}$ of $F_i,E_i$ such that $F_i^a = (F_i^{(a)})^{\oplus [a]!}$ and $E_i^{a} = (E_i^{(a)})^{\oplus [a]!}$.  
The categorical action is supplied by
\cite[Theorems~5.2 and~6.2]{KK}. 
For \(H=\langle \alpha^\vee_i,\mu\rangle >0\), the functor 
\[
 F_i^{(H)}:\mathcal C_\mu\longrightarrow\mathcal C_{s_i\mu}
\]
is an equivalence if $\mathcal{C}_{\mu +\alpha_i} =0$ and we call it a \emph{Scopes equivalence}.
In this case, a quasi-inverse is given by certain shift of $E_i^{(H)}$.
We refer to \cite[Section 3]{Web23} for more details.

\subsection{The exact abacus criterion}
Let \(\bar\alpha_i\in \mathring{\mathfrak{h}}^*_{\mathbb{R}}\) be the finite part of \(\alpha_i\), so
\[
 \bar\alpha_0=-2\eps_1,\qquad \bar\alpha_i=\eps_i-\eps_{i+1}\ (0<i<e),\qquad
 \bar\alpha_e=2\eps_e.
\]
Recall that \(\mu=\ell\Lambda_0+P-c_0\delta\), so 
\begin{equation}\label{eq:coroot}
 H_i = \langle \alpha^\vee_i,\mu \rangle =\begin{cases}
     \ell - P_1 & \text{if } i =0; \\
     P_i - P_{i+1} & \text{if } 0 < i < e; \\
     P_e & \text{if } i=e.
 \end{cases}
\end{equation}

\begin{theorem}\label{thm:scopes}
Let \(\mu \in P(\Lambda)\) be with \(H_i> 0\). 
Then $F_i^{(H_i)}|_{\mathcal{C}_\mu}$ is a Scopes equivalence if and only if
\begin{equation}\label{eq:scopes}
 c_0-\delta_{i,0}<T_{\mathbf k}(P+ \bar\alpha_i).
\end{equation}
Equivalently, every multipartition of content \(\Lambda -\mu\) has no
removable \(i\)-nodes.
\end{theorem}
\begin{proof}
Set $\nu = \mu + \alpha_i = \ell \Lambda_0 +P' -c'_0$, then $P'=P+\bar\alpha_i$ and $c'_0 =c_0-\delta_{i,0}$.
By \cref{prop:support}, $\nu\in P(\Lambda)$ if and only if $\exists \blam$ such that $\Lambda -\nu =\cont_{\mathbf{k}}(\blam)$. 
The latter holds if and only if $c'_0\ge T_{\mathbf{k}}(P')$ by \cref{thm:higher} and reconstruction algorithm.
This completes the proof of \eqref{eq:scopes}.

The node formulation then follows from the $U_q(\mathfrak{g})$-action on the Fock space $\mathcal{F}(\Lambda)$, see \cite[Theorem~6D.3(b)]{EM24}.
\end{proof}

\begin{corollary}\label{cor:scopes-level-one}
Assume that $\Lambda = \Lambda_k$, let $p\in \mathcal{P}_k$ and $\mu = \Lambda_0 +p-(f_k(p)+m)\delta$. 
If $H = \langle \alpha^\vee_i,\mu\rangle \neq 0$ for $i\in \{ 0,\ldots ,e\}$.  
Let \(h=h(p)\) and \(h'=h(p+\epsilon\bar\alpha_i)\) with $\epsilon = \operatorname{sgn}(H)$. 
Then $\mathcal{C}_{\mu}$ and $\mathcal{C}_{s_i(\mu)}$ are Scopes equivalent if and only if
\begin{equation}\label{eq:scopes-one}
 |H| \ge \begin{cases}
     m & \text{if } i=0,e,\\
     2m + \dfrac{|h-k|-|h'-k|}{2} & \text{otherwise.}
 \end{cases}
\end{equation}
\end{corollary}
\begin{proof}
We may assume that $H>0$ since $\langle \alpha^\vee_i,s_i(\mu)\rangle =-H$. 
At level one, \eqref{eq:scopes} is equivalent to $m<f_k(p + \bar\alpha_i)-f_k(p)+ \delta_{i,0}$ for some/any $\lambda$ of weight $\Lambda_0+P-c_0\delta$. By definition,
\[
f_k(p+\bar\alpha_i) - f_k(p) + \delta_{i,0} = \begin{cases}
    2-p_1 & \text{if } i =0, \\
    \frac{2(p_i - p_{i+1}+1)+|h'-k|-|h-k|}{4} & \text{if } 0 < i < e, \\
    p_e+1 & \text{if } i=e.
\end{cases}
\]
The statement then follows from \eqref{eq:coroot}.
\end{proof}


\section{Scopes normal forms and counting}\label{sec:counting}

We now work at level one. Fix \(0\le k,t\le e\) with $\Lambda =\Lambda_k$ and
\(t\equiv k\pmod2\), and \(m\ge0\). The weights with
\(h(p)=t\) and depth \(m\) form one $W$-orbit.
Recall from the last section that Scopes equivalences induce equivalence relation on $W$-orbits in $P(\Lambda)$.
We call such an equivalence class Scopes class.

For \(a,b\in\mathbb{Z}\), define
\begin{equation}\label{eq:threshold}
 L_{ab}=2m+
 \begin{cases}
 2-a'-b',&t<k,\\
 0,&t=k,\\
 a'+b',&t>k.
 \end{cases}
\end{equation}
where $a',b'\in \{0,1\}$ with $a\equiv a'\pmod 2$ and $b\equiv b'\pmod 2$.
The dependence on \(k,t,m\) will always be understood.

\begin{theorem}\label{thm:normal-forms}
Suppose \(m\ge1\). Every Scopes class contains exactly one weight $\mu$ with vector
\(p\) satisfying
\begin{equation}\label{eq:terminal}
 p_1\le m,\qquad p_e\ge1-m,\qquad
 p_{i+1}-p_i< L_{p_i,p_{i+1}}\quad(1\le i<e).
\end{equation}
When \(m=0\), there is one Scopes class, with terminal vector \(u(t)\).
\end{theorem}
\begin{proof}
We prove the theorem by showing that the desired $\mu$ is the unique weight in this class that has minimal content.
Start with arbitrary weight $\mu^{(0)} =\Lambda_0 + p^{(0)}-c^{(0)}_0\delta$ in the fixed Scope class, if there is $i\in \{ 0,\ldots ,e\}$ such that $\langle h_i,\mu^{(0)}\rangle <0$ and $\mu^{(0)}-\alpha_i\not\in P(\Lambda)$, then $\mu^{(1)}=s_i(\mu^{(0)})$ also belongs to this Scope class and $\operatorname{ht}(\Lambda_k -\mu^{(0)})>\operatorname{ht}(\Lambda_k-\mu^{(1)})$. 
Repeating this step, we obtain \(\mu^{(0)},\mu^{(1)},\mu^{(2)},\ldots\).
This process must terminate because the height strictly decreases at each step.
We denote the terminal weight by \(\mu'=\Lambda_0+p'-(f_k(p')+m)\delta\).

For $i\in \{0,\ldots ,e\}$, if $H'_i =\langle \alpha^\vee_i ,\mu'\rangle <0$, then $s_i(\mu') + \alpha_i\in P(\Lambda_k)$ by the minimality of $\mu'$. 
By \cref{cor:scopes-level-one}, we have
\[
\begin{cases}
    1-p'_1 > -m & \text{if } i=0, \\
    p'_i-p'_{i+1} > -2m+\frac{|t+(-1)^{p'_i}+(-1)^{p'_{i+1}}-k|-|t-k|}{2} & \text{if } 0<i<e,\\
    p'_e > -m & \text{if } i=e.
\end{cases}
\]
This is obviously equivalent to \eqref{eq:terminal}. 

It leaves to show any two moves $(\mu^{(0)},\mu^{(1)}\ldots)$ and $(\nu^{(0)},\nu^{(1)},\ldots)$ with $\mu^{(0)} = \nu^{(0)}=\nu$ terminate at the same weight.
Let $\ell (\nu)$ be the smallest number of steps among such moves stated with $\nu$.
We prove the uniqueness by induction on $\ell (\nu)$. 
The statement holds trivially when $\ell (\nu) =0$.
Now, assume that $n>0$ and the uniqueness has been proven for all $\nu'$ with $\ell (\nu')<n$. 
Let $i,j\in \{ 0,\ldots ,e \}$ s.t. $\mu^{(1)} = s_i(\mu^{(0)})$ and $\nu^{(1)} =s_j(\nu^{(0)})$. 

We claim that for each $w\in W' :=\langle s_i,s_j\rangle$, $w(\nu)$ belongs to the given Scopes class.
Once this claim is proved, we take a move $(\widetilde{\nu}^{(0)},\widetilde{\nu}^{(1)}\ldots )$ started with $w_0(\nu)$ for the longest element $w_0\in W'$, then
\[
(\nu, s_i(\nu), s_js_i(\nu),\ldots ,w_0(\nu),\widetilde{\nu}^{(1)},\widetilde{\nu}^{(2)}\ldots)
\]
and 
\[
(\nu, s_j(\nu), s_is_j(\nu),\ldots ,w_0(\nu),\widetilde{\nu}^{(1)},\widetilde{\nu}^{(2)}\ldots)
\]
are moves started with $\nu$. 
Applying the induction hypothesis to $s_i(\nu)$ and $s_j(\nu)$ respectively, we see that the two given moves terminate at the same weight.
Without loss of generality, assume that $i\leq j$.
Now we turn to prove the claim.

1) $i\neq j- 1$: It is clear.

2) $0<i<i+1=j<e$: $p_{i+2}-p_{i+1}\ge L_{p_{i+1},p_{i+2}}$ and $p_{i+1}-p_{i+1}\ge L_{p_i,p_{i+1}}$, then $s_i(\nu)\sim s_{i+1}s_i(\nu)\sim s_{i}s_{i+1}s_i(\nu)$ can be verified by applying the following inequalities to \cref{cor:scopes-level-one}: 
\begin{align*}
    & p_{i+2} - p_i \ge L_{p_{i},p_{i+1}} + L_{p_{i+1},p_{i+2}} > L_{p_{i},p_{i+2}} \quad (\Leftarrow m\ge 1), \\
    & p_{i+2} - p_{i+1} \ge L_{p_{i+1},p_{i+2}} ,
\end{align*}
The other chain $s_{i+1}(\nu)\sim s_is_{i+1}(\nu)\sim s_{i+1}s_is_{i+1}(\nu)$ can be proved similarly.

3) $i=0,j=1$: $p_1 > m\ge 1$ and $p_2 - p_1 \ge L_{p_{1},p_{2}}$. The chain $s_0(\nu)\sim s_1s_0(\nu) \sim s_0s_1s_0(\nu) \sim w_0(\nu)$ then follows from
\begin{align*}
    & p_2-(2-p_1) > p_2 - p_1 \ge L_{2-p_1,p_2} \quad (\Leftarrow p_1> 1) \\
    & p_2 \ge p_1 + L_{p_1,p_2} \ge p_1 > m \\
    & (2-p_1) - (2-p_2) = p_2 - p_1 \ge L_{2-p_2,2-p_1}
\end{align*}
Similar for the other chain.

4) $i=e-1,j=e$: similar to 3) and we leave it to the readers.

The case $m =0$ is a direct consequence of \cref{cor:scopes-level-one}.
\end{proof}

The later product formula gives the cardinality of this set of normal forms.

\begin{theorem}\label{thm:count}
For \(m\ge1\), the number of Scopes classes in the orbit
\(h(p)=t\), \(D_k=m\), is
\begin{equation}\label{eq:scopes-count}
 \begin{cases}
 \displaystyle\binom{em}{t}\binom{em+e-t-1}{e-t},&t<k,\\[6pt]
 \displaystyle\binom{em}{t}\binom{em}{e-t},&t=k,\\[6pt]
 \displaystyle\binom{em+t-1}{t}\binom{em}{e-t},&t>k.
 \end{cases}
\end{equation}
For \(m=0\) the number is one.
\end{theorem}
\begin{proof}
At positive depth apply \cref{thm:level-one-factorization}
and divide its chamber count by $2^e t!(e-t)!$, as in
\cref{thm:all-chambers}. The displayed products simplify to
\eqref{eq:scopes-count}. The zero-depth assertion is
\cref{thm:normal-forms}.
\end{proof}

\begin{remark}\label{rem:catalan}
When \(k=t=0\), the count is \(\binom{em}{e}\). In the usual
normalization
\[
 \operatorname{Cat}^{r}(C_e)
   =\prod_{i=1}^e\frac{2er+2i}{2i}
   =\binom{(r+1)e}{e},
\]
this is \(\operatorname{Cat}^{m-1}(C_e)\); the parameter is
\(m-1\), not \(m\). The case \(k=t=e\) is the same.
Fishel--Vazirani relate dominant Shi regions to simultaneous cores
in type A, and Armstrong--Hanusa--Jones give the corresponding
type-C construction with self-conjugate cores \cite{FV,AHJ}.
Our normal forms give concrete class representatives; no identification
with those simultaneous-core sets is asserted here.
\end{remark}

\begin{example}\label{ex:four-classes}
For $e=2$ and $k=t=m=1$, the four Scopes classes have terminal
vectors $(-1,0)$, $(0,1)$, $(1,0)$ and $(1,2)$.
Their RoCK classes are identified in \cref{ex:rouquier}.
\end{example}

\section{Scopes arrangements and counting}\label{sec:arrangements}

\subsection{The arrangement from normalized coordinates}\label{sec:ArrangementAndNormalized}
Return to $\Lambda$ of arbitrary level \(\ell\). Let
\[
 \mathring{\Phi}^+=\{2\eps_i:1\le i\le e\}\cup
 \{\eps_i-\eps_j,\eps_i+\eps_j:1\le i<j\le e\},
 \qquad \mathring{\Phi}=\mathring{\Phi}^+\cup(-\mathring{\Phi}^+).
\]
be the positive root system and root system of type $C_e$.
Then positive real roots in $\Phi_{\mathrm{re}}(C_e^{(1)})$ are of the form \(\gamma+n\delta\), with
\(n\ge0\) for \(\gamma\in\mathring{\Phi}^+\), and \(n\ge1\) otherwise.
On the level-\(\ell\) slice their hyperplanes are
\[
 H_{\gamma+n\delta}=\{x\in \mathring{\mathfrak{h}}^*_{\mathbb{R}}:(\gamma,x)+2\ell n=0\}.
\]
Let \(\mathfrak a_\ell=\{\ell>x_1>\cdots>x_e>0\}\). For \(\mu \in P^+(\Lambda)\), put
\[
 N_\mu=\{\alpha\in \Phi^+_{\mathrm{re}}:\mu+\alpha\in P(\Lambda)\}.
\]
Webster \cite{Web23} shows that this is a finite hyperplane arrangement of $\mathring{\mathfrak{h}}^*_{\mathbb{R}}$. 
The connected components of \(\mathring{\mathfrak{h}}^*_{\mathbb{R}}\setminus\bigcup_{\alpha\in N_\mu}H_\alpha\)
are called \emph{Scopes chambers}. 
A chamber is \emph{RoCK} if it contains a translate of a finite Weyl chamber. 
There is a bijection between Scopes classes in $W\mu$ and $W_\mu$-orbits of Scopes chambers with respect to $N_\mu$: the weight $\nu = w^{-1}(\mu)$ corresponds to the $W_\mu$-orbit of the chamber containing $w\mathfrak{a}_\ell$.
We call the category $\mathcal{C}_\nu$ RoCK if the corresponding alcove lies in a RoCK chamber.
This definition applies to any integrable categorical action with support \(P(\Lambda)\).

\begin{proposition}\label{prop:walls}
Write \(\mu=\ell\Lambda_0+a-C\delta\), where
\(\ell\ge a_1\ge\cdots\ge a_e\ge0\), and \(C=\zeta_\Lambda(a)+M\). 
For \(\gamma\in\mathring{\Phi}\), put
\[
 b_\gamma=C-T_{\mathbf k}(a+\gamma).
\]
Then \(N_\mu\) consists exactly of
\begin{equation}\label{eq:wall-cutoffs}
 \gamma+n\delta\quad\text{with}\quad
 \begin{cases}0\le n\le b_\gamma,&\gamma\in\mathring{\Phi}^+,\\
 1\le n\le b_\gamma,&\gamma\in-\mathring{\Phi}^+.
 \end{cases}
\end{equation}
An interval with upper endpoint below its lower endpoint is empty.
Thus \(2e^2\) evaluations of \eqref{eq:closed-T} determine the entire arrangement, in finite interval form.
\end{proposition}
\begin{proof}
The weight \(\nu =\mu+\gamma+n\delta\) has coordinates
\((a+\gamma,C-n)\). \cref{thm:higher} says exactly that
$\nu\in P(\Lambda)$ when \(C-n\ge T_{\mathbf k}(a+\gamma)\).
\end{proof}

The integers $b_\gamma$ in our notation are denoted by $k_i^\pm$ in Webster's paper \cite{Web23}.

\begin{proposition}\label{prop:wall-growth}
Let $a$ be that of \cref{prop:walls}. 
For \(\gamma\in\mathring{\Phi}^+\), write
\[
 d_\gamma^\pm=T_{\mathbf k}(a\pm\gamma)-\zeta_\Lambda(a),
 \qquad \kappa_\gamma=1-d_\gamma^+-d_\gamma^-.
\]
Then \(0\le d_\gamma^+\le2\), \(-1\le d_\gamma^-\le1\), and
\(\kappa_\gamma\in\{-1,0,1,2\}\). 
\end{proposition}
\begin{proof}
Put \(\nu=\ell\Lambda_0+a-\zeta_\Lambda(a)\delta\).
This is a dominant maximal weight in $P(\Lambda)$.
Note that $b_{\pm \gamma}(M) =M-d_{\gamma}^{\pm}$.
Applying \cref{prop:walls} to $M\gg 0$, we know that $d_\gamma^{\pm}$ is the smallest integer such that
\[
\nu \pm \gamma - d^{\pm}_{\gamma}\delta \in P(\Lambda) .
\]

Recall that 
\begin{equation}\label{eq:lowering_+realroot}
    \mu \in P(\Lambda), \alpha \in\Phi_{\mathrm{re}}^+ , (\mu, \alpha) >0 \Rightarrow \mu -\alpha \in P(\Lambda) 
\end{equation}

The upper bounds of $d_\gamma^{\pm}$ then follow immediately.
Suppose that $d^+_{\gamma} \le -1$, then $\mu = \nu + \gamma + \delta\in P(\Lambda)$. Since $(\mu,\gamma) \ge (\gamma,\gamma)\ge 2$, we have $\nu + \delta\in P(\Lambda)$ by \eqref{eq:lowering_+realroot}, which contradicts to the maximality of $\nu$. This proves the lower bound of $d_\gamma^+$. The proof for $d_\gamma^-$ is the same.


If \(d_\gamma^+=2\), lowering \(\nu\) by
\(\delta-\gamma\) is impossible. 
So $(\nu,\delta -\gamma)= 0$ and thus $\langle \gamma^\vee , \nu \rangle >0$.
Therefore, \(\nu-\gamma \in P(\Lambda)\) and
\(d_\gamma^-\le0\). This excludes the only possibility for
\(\kappa_\gamma<-1\); the other bounds give
\(\kappa_\gamma\le2\). 
\end{proof}

\begin{lemma}\label{lem:window}
Write \(\beta =\Lambda-\nu=\sum_{i=0}^e c_i\alpha_i\), where
\(\nu=\ell\Lambda_0+a-\zeta_\Lambda(a)\delta\). Then
\[
 \begin{cases}
 \kappa_\gamma\in\{-1,0,1\}\quad\text{for every }\gamma\in\mathring{\Phi}^+,
       & \text{if }\min_i c_i=0,\\
 \kappa_\gamma\in\{0,1,2\}\quad\text{for every }\gamma\in\mathring{\Phi}^+,
       & \text{if } \min_ic_i > 0.
 \end{cases}
\]
\end{lemma}
\begin{proof}
If \(\kappa_\gamma=2\), the bounds already proved force \((d_\gamma^+,d_\gamma^-)=(0,-1)\). 
Hence \(\nu+\gamma,\nu+\delta-\gamma \in P(\Lambda)\) 
and \(2\beta-\delta\in Q_+\backslash \{ 0\}\). 
Since the coefficients of \(\delta\) are \(1,2,\ldots,2,1\), every \(c_i\) is positive.

Suppose now that \(\min_ic_i>0\). 
At level one, write \(a=(1^t,0^{e-t})\) and \(\Lambda=\Lambda_k\).
Then \(c_0=\max(t-k,0)/2\) and \(c_e=\max(k-t,0)/2\), 
so our assumption implies \(\ell\ge2\).
We claim that \(\nu+\alpha \in P(\Lambda)\) for $\alpha\in \Phi_{\mathrm{re}}$ whenever \(s_{\alpha}\in W_\nu\).
We only need to prove for $\alpha >0$ since $\nu - \alpha = s_\alpha (\nu + \alpha)$ and $P(\Lambda)$ is $W$-stable. 

Note that $W_\nu= \langle s_i \mid I_\nu \rangle \subsetneq W$, where $I_\nu = \{ i\in I \mid  a_i = a_{i+1}\}$, where we use the convention $a_0 = \ell$ and $a_{e+1} =0$.
After applying the diagram automorphism $i \longleftrightarrow e-i$, we may assume that $\alpha \in \mathring{\Phi}^+$. We have three cases to deal with:

(1) $\alpha = \varepsilon_i - \varepsilon_j$ with $i< j$ and $a_i = a_j > 0$. 
The dominant weight in $W(\nu + \alpha)$ is $\mu =r_{i',i}r_{j,j'}(\nu + \alpha)=\nu + \varepsilon_{i'} - \varepsilon_{j'}$, where $i' = \min \{t \mid a_t =a_i \}$, $j'=\max \{ t \mid a_t = a_j\}$ and $r_{a,b}$ is the permutation $(a,b)$. 
Since $\varepsilon_{i'} - \varepsilon_{j'}=\sum_{t=i'}^{j'-1}\alpha_t$, $\Lambda - \mu\in Q_+$, and hence $\nu + \alpha \in P(\lambda)$.

(2) $\alpha = \varepsilon_i + \epsilon \varepsilon_j$ with $\epsilon \in \{\pm 1\}$, $i< j$ and $a_i = a_j =0$, then the dominant weight in $W(\nu + \alpha)$ is $\mu = \nu + \varepsilon_{i'} + \varepsilon_{i'+1}$ with $i' = \min \{ t \mid a_t =0\}$. 
Since $\varepsilon_{i'}+\varepsilon_{i'+1} =\alpha_{i'}+2\sum_{t = i'+1}^{e-1}\alpha_t + \alpha_e$, we only need to prove that $c_t \ge 2$ for those $t<e$ with $a_t =0(t< e)$.
Indeed, $\sum_{s}a_{ts}c_s =\langle \alpha_t^{\vee} , \beta \rangle =\langle \alpha_t^\vee , \Lambda\rangle \ge 0$. 
The statement $c_t \ge 2$ then follows from a downward induction argument in $t$.

(3) $\alpha = 2\varepsilon_i$ with $a_i =0$, then the dominant weight in $W(\nu + \alpha)$ is $\nu + \varepsilon_{i'}+\varepsilon_{j'}$ with $i' = \min \{ t \mid a_t \le 1 \}$ and $j'=\min \{t \mid a_t =0\}$. 
The rest is the same as (2).

Finally, \(\kappa_\gamma=-1\) forces
\((d_\gamma^+,d_\gamma^-)=(1,1)\) or \((2,0)\).
In the first case, $\nu - \gamma \not\in P(\Lambda)$ and hence $(\gamma ,\nu) = 0$ by \eqref{eq:lowering_+realroot}.
Similarly in the second case, $(\nu , \delta - \gamma) =0$ and $\nu + \delta -\gamma\not\in P(\Lambda)$.
The claim forces $\min_i c_i =0$.
\end{proof}

In the rest of this section, we fix $\Lambda$ and $a$ with $\ell \ge a_1 \ge \ldots \ge a_e\ge 0$. 
Let \(\mathcal A_{a,M}\) be the affine arrangement determined by $N_{\nu}$ with $\nu = \ell\Lambda_0 +a -(\zeta_{\Lambda}(a)+M)\delta$. 
Its characteristic polynomial is
\[
 \chi_{a,M}(q)=\sum_{X\in L(\mathcal A_{a,M})}
       \mu_L(\mathring{\mathfrak{h}}^*_{\mathbb{R}},X)q^{\dim X},
\]
where \(L(\mathcal A_{a,M}) = \{ \left(\bigcap_{H\in\mathcal H} H\right)\cap \mathring{\mathfrak{h}}^*_{\mathbb{R}}\neq \emptyset \mid \mathcal{H}\subset \mathcal{A}_{a,M} \}\),
including \(E\), ordered by reverse inclusion, and \(\mu_L\)
is its M\"obius function.

\begin{theorem}\label{thm:all-chambers}
Let \(V(\Lambda)\) be an integrable highest-weight module of an untwisted affine Lie algebra of finite rank \(r\), at positive level.
Let \(\eta\in P(\Lambda)\) be dominant, and suppose that
\(\eta+\delta\in P(\Lambda)\). 
The finite parabolic group \(W_\eta\) acts freely on the set of Scopes chambers with respect to Scopes arrangement
\[
 \mathcal A_\eta=
 \{H_\beta:\beta\in \Phi_{\mathrm{re}}^+,\ \eta+\beta\in  P(\Lambda)\}.
\]
The number of
Scopes classes in \(W\cdot\eta\) is therefore
\begin{equation}\label{eq:general-chambers}
 |\operatorname{Scopes}(W\eta)|
       =\frac{(-1)^r\chi_{\mathcal A_\eta}(-1)}{|W_\eta|}.
\end{equation}
\end{theorem}
\begin{proof}
Choose a positive real root \(\beta\) with \(s_\beta\in W_\eta\).
Thus \((\eta,\beta)=0\). Untwistedness implies that
\(\delta-\beta\) is real, and
\[
 \langle\eta+\delta,(\delta-\beta)^\vee\rangle
       =\frac{2(\eta,\delta)}{(\beta,\beta)}>0.
\]
Then $\eta + \beta = \eta + \delta -(\delta -\beta)\in P(\Lambda)$ by \eqref{eq:lowering_+realroot}. 
Hence \(\beta\in N_\eta\).

For $\alpha \in N_\eta$ and $\beta\in \Phi_{\mathrm{re}}^+$ with $s_\beta \in W_\eta$, if $s_\beta (\alpha) <0$, then $(\eta ,s_\beta (\alpha)) = (\eta , \alpha) =0$. 
Hence, $\mathcal{A}_\eta$ is $W_\eta$-invariant and so is the set of Scopes chambers. 
Combining this with the above, we obtain that $W_\eta$ acts freely on the set of Scopes chambers.

All reflecting hyperplanes of \(W_\eta\) occur among the Scopes
walls. Every Scopes chamber consequently lies in a unique open
\(W_\eta\)-chamber. Since a finite reflection group acts freely
on its open chambers, it acts freely on the Scopes chambers.
The identity \eqref{eq:general-chambers} then follows from the Zaslavsky's theorem \cite[Theorem A]{Zaslavsky}.
\end{proof}

\begin{corollary}\label{cor:typeC-chambers}
In type \(C_e^{(1)}\), at depth \(M\ge1\), set
\(n_j=\#\{i:a_i=j\}\). Then
\begin{equation}\label{eq:all-chambers}
 |\operatorname{Scopes}(W\mu)|
 =\frac{(-1)^e\chi_M(-1)}
        {2^{n_0+n_\ell}\prod_{j=0}^{\ell}n_j!}.
\end{equation}
The same formula holds at \(M=0\) if every simple-root coefficient
of \(\Lambda-\nu\) is positive.
\end{corollary}
\begin{proof}
The stabilizer $W_\mu$ is the standard parabolic subgroup of type
\(C_{n_0}\times C_{n_\ell}\times
\prod_{j=1}^{\ell-1}A_{n_j-1}\) whose order is $2^{n_0+n_\ell}\prod_{j=0}^\ell n_j!$ as desired.
Its order gives the denominator. 
When $D_{\mathbf{k}}(\mu) >0$, \(\mu+\delta \in P(\Lambda)\), so \cref{thm:all-chambers} applies.

In the stated depth-zero case, the proof of
\cref{lem:window} implies that $\beta \in N_\mu$ if $s_\beta\in W_\mu$, so $W_\mu$ acts freely on the set of Scopes chambers.
The proof of \eqref{eq:general-chambers} applies in this case.
\end{proof}


\begin{proposition}[{\cite[Theorem 2.2]{Athanasiadis}}]\label{prop:finite-field}
Let \(\mathcal A\) be a rational affine arrangement in dimension \(r\),
and choose integer equations \(v_H\cdot x=b_H\) for its hyperplanes.
Exclude the primes dividing any nonzero square minor of the augmented
matrix with rows \((v_H,b_H)\). For every remaining prime \(p\) and
every power \(q=p^s\),
\[
 \chi_{\mathcal A}(q)=
 |\mathbb F_q^r\setminus\textstyle\bigcup_{H\in\mathcal A}H(\mathbb F_q)|.
\]
Thus finitely many finite-field counts determine the exact Scopes-class
count in \eqref{eq:general-chambers}.
\end{proposition}

\subsection{Band decomposition}\label{sec:bands}
In this subsection let \(V(\Lambda)\) be an integrable irreducible highest-weight module of positive level for an untwisted affine algebra.
Let \(\mathring{\Phi}\) be the underlying finite root system of rank $r$, \(\theta =\sum_{i=1}^r c_i\alpha_i\) the highest root of $\mathring{\Phi}$, and \(\delta\) the null root. 
Recall that for \(\eta \in P^+(\Lambda)\), 
\[
 N_\eta=\{\alpha\in \Phi_{\mathrm{re}}^+:\eta+\alpha\in P(\Lambda)\}.
\]
Fix \(\nu\in P^+(\Lambda)\). 
Define the finite sets
\[
 \begin{split}
 & O_\nu=\{\alpha\in N_\nu:(\nu,\alpha)=0\},\\
 & B_M(\nu)=\{\alpha \in \Phi_{\mathrm{re}}^+:(\nu,\alpha -M\delta)<0\}
 \qquad(M\in\Z_{\ge1}).
 \end{split}
\]
Since $\alpha\in\Phi^+_{\mathrm{re}}$ is of the form $\beta + n\delta$ with $\beta\in\mathring{\Phi}$ and $n\ge 0$, the set $B_M(\nu)$ is finite.

\begin{theorem}\label{thm:band-decomposition}
For every dominant supported \(\nu\) and every integer \(M\ge1\),
\begin{equation}\label{eq:band-decomposition}
 N_{\nu-M\delta}=B_M(\nu)\ \sqcup\
       (M\delta-O_\nu)\ \sqcup\ (M\delta+N_\nu).
\end{equation}
\end{theorem}

\begin{proof}
It is obvious that the union is  disjoint and $M\delta + N_\nu \subset N_{\nu - M\delta}$.
If $\beta \in B_M(\nu)$, then $\alpha \in N_{\nu - M\delta}$ by applying  \eqref{eq:lowering_+realroot} to $\mu = \mu =\nu$ and $\alpha = M\delta -\beta$.
For $\beta \in O_\nu$, $\nu - \beta = s_\beta (\nu + \beta)\in P(\Lambda)$, so $M\delta -\beta \in N_{\nu - M\delta}$.
This completes the proof of $\supseteq$.


In contrast, let $\alpha \in N_{\nu - M\delta}$ and $\beta =\alpha -M\delta$. 
If $(\nu , \beta) < 0$, then $\alpha \in B_{M}(\nu)$.
If $\beta > 0$, then $\nu + \beta \in P(\Lambda)$ and $\beta \in N_{\nu}$, thus $\alpha \in N_{\nu} + M\delta$.
If $\beta < 0$ and $(\nu , \beta) = 0$, then $\nu - \beta = s_\beta (\nu + \beta) \in P(\Lambda)$, hence $-\beta \in O_\nu$ and $\alpha \in M\delta - O_\nu$. This proves $\subseteq$.
\end{proof}

Write $\nu = \ell \Lambda_0 + K\bar{p}_{\nu} -c\delta$ where $\bar{p}_\nu\in\mathring{\mathfrak{h}}^*$ and \( K =(\nu,\delta) >0 \).
The vector \(\bar{p}_\nu\) lies in the closed fundamental alcove, in particular, $0\le (\alpha_i ,\bar{p}_\nu) \le (\theta ,\bar{p}_\nu)\le 1$ for $i\neq 0$. 
Then $B_M(\nu)$ gives the following set of hyperplanes in $\mathring{\mathfrak{h}}^*_{\mathbb{R}}$:  
\begin{equation}\label{eq:open-band}
 \mathcal B_M(\bar{p}_\nu) = \{ H_{\alpha,j}:\alpha\in\mathring{\Phi}^+,\ j\in\Z, \ |j-(\alpha,\bar{p}_\nu)|<M \}
\end{equation}
where $H_{\alpha,j} = H_{-\alpha +j\delta} = \{z\in\mathring{\mathfrak{h}}^*_{\mathbb{R}}: (\alpha ,z)=j\}$.
For \(\nu=\Lambda\), the sets \(N_\Lambda\) and
\(O_\Lambda\) are empty.
Thus the Scopes arrangement of \(\Lambda-M\delta\) is precisely \(\mathcal B_M(\bar{p}_\Lambda)\).
The arrangement depends only on the relatively open alcove face
containing \(p_\Lambda\).


\subsection{Level-one vertex bands in type C}
Return to type \(C_e^{(1)}\) and fix \(\Lambda=\Lambda_k\).
For \(0\le t\le e\) with \(t\equiv k\pmod2\), put
\[
 a_t=(1^t,0^{e-t}),\qquad
 \zeta_t=\frac{\max(t-k,0)}2,\qquad
 \nu_t=\Lambda_0+a_t-\zeta_t\delta,\qquad p_t=\frac{a_t}{2}.
\]
Thus, \(\nu_t\) is the dominant maximal weight corresponding to the alcove vertex \(a_t\). 
Define
\[
 E_{k,t}^+=
 \begin{cases}
 \{\eps_i-\eps_j,\eps_i+\eps_j:1\le i<j\le t\},&t>k,\\
 \{\eps_i-\eps_j,\eps_i+\eps_j:t<i<j\le e\},&t<k,\\
 \varnothing,&t=k.
 \end{cases}
\]
The nonempty sets are the positive short roots of the indicated
\(C_t\) or \(C_{e-t}\) subsystem.

\begin{proposition}\label{prop:vertex-bands}
For \(M\ge1\), the Scopes walls of \(\nu_t-M\delta\) of $\mathcal{A}_{\nu_t - M\delta}$ are
exactly \(H_{\gamma,j}\), with
\(\gamma\in\mathring{\Phi}^+\), \(j\in\Z\) satisfying
\begin{equation}\label{eq:vertex-bands}
 \begin{cases}
 |j-(\gamma,p_t)|\le M,&\gamma\in E_{k,t}^+,\\
 |j-(\gamma,p_t)|<M,&\gamma\notin E_{k,t}^+.
 \end{cases}
\end{equation}
\end{proposition}
\begin{proof}
We compute the thresholds
\(d_\gamma^\pm=T_{(k)}(a_t\pm\gamma)-\zeta_t = f_k(a_t\pm \gamma) - \zeta_t\).
By \cref{thm:depth}, a direct computation gives the following list 
\begin{equation*}
\begin{array}{l|l|l}
\toprule
 \gamma&\text{condition}&(d_\gamma^+,d_\gamma^-)\\
\midrule
 2\eps_i&i\le t&(2,0)\\
 2\eps_i&i>t&(1,1)\\
 \eps_i\pm\eps_j&i\le t<j&(1,0)\\
 \eps_i-\eps_j&i<j\le t,\ t>k&(0,0)\\
 \eps_i-\eps_j&i<j\le t,\ t\le k&(1,1)\\
 \eps_i+\eps_j&i<j\le t,\ t>k&(1,-1)\\
 \eps_i+\eps_j&i<j\le t,\ t\le k&(2,0)\\
 \eps_i\pm\eps_j&t<i<j,\ t<k&(0,0)\\
 \eps_i\pm\eps_j&t<i<j,\ t\ge k&(1,1)\\
\bottomrule
\end{array}
\end{equation*}
Combining with \cref{prop:walls}, one can check that $\alpha \in N_{\nu_t}$ must be of the form $\pm (\gamma - (\gamma ,p_t)\delta)$ for some $\gamma \in E_{k,t}^+$. 
Note that $(\gamma ,p_t) \in \{ 0,1\}$ for $\gamma\in E^+_{k,t}$ and $(\gamma ,p_t) = 1$ only if $\gamma = \varepsilon_i + \varepsilon_j$ with $i<j\le t$.
So, $M\delta + N_{\nu_t}\bigsqcup M\delta - O_{\nu_t}$ contributes to walls of the form $H_{\gamma ,j}$ with $\gamma \in E_{k,t}^+$ and $j - (\gamma ,p_t)=\pm M$. 
The result then follows from \eqref{eq:band-decomposition} and \eqref{eq:open-band}. 
\end{proof}

For reference, the number of walls in this arrangement is
\[
 2Me^2-t^2-(e-t)^2+
 \begin{cases}
 2t(t-1),&t>k,\\
 2(e-t)(e-t-1),&t<k,\\
 0,&t=k.
 \end{cases}
\]
The first three terms count the open band; the remaining term counts
the two endpoints for each selected short-root direction. 

\subsection{Highest-weight factorization in every untwisted type}
\label{sec:affine-band-factorization}
Let us go back to the general untwisted setting as in \cref{sec:bands}.
Write $\hcox=1+\sum_i c_i$ for its Coxeter number and
\[
 \mathfrak a=\{x\in \mathring{\mathfrak{h}}_{\mathbb{R}}:\alpha_i(x)>0\ (1\leq i\leq e),\ \theta(x)<1\}
\]
for the \emph{open} fundamental alcove; its closure
$\overline{\mathfrak a}$ is defined by the corresponding non-strict
inequalities. 
For $p\in\overline{\mathfrak a}$ and an integer
$M\geq1$, we define
\begin{equation}\label{eq:band-definition}
 \mathcal B_M(p)=
 \bigl\{H_{\alpha ,j}:\alpha\in\mathring{\Phi}^+,\ j\in\mathbb Z,
                  \ |j-\alpha(p)|<M\bigr\}.
\end{equation}
Let $W_F$ be the stabilizer of $p$ in the affine Weyl group, where $F$
is the facet of $\overline{\mathfrak a}$ containing $p$ in its interior.
It is a finite parabolic subgroup.  
Denote its essential rank by $r_F$ and its exponents, with multiplicities, by $e_1,\ldots,e_{r_F}$.

\begin{theorem}\label{thm:affine-band-factorization}
For every $p$, and positive integer $M$ as above,
\begin{equation}\label{eq:band-characteristic-polynomial}
 \chi(\mathcal B_M(p),q)
   =(q-M\hcox)^{r-r_F}\prod_{i=1}^{r_F}(q-M\hcox+e_i).
\end{equation}
The empty product is $1$.  
In particular, the formula applies to every face, including non-special vertices, in every irreducible
crystallographic type.
\end{theorem}

The proof is a face-centred version of the classical alcove count:
compare Haiman \cite[\S7.4]{Haiman}, Athanasiadis
\cite[\S3]{AthanasiadisCatalan} and Sommers \cite[Theorem~5.7]{Sommers}.
The stabilizer calculation below keeps the centre on an arbitrary facet. 
Write $\Pcoweight$ and $\Qcoroot$
for the coweight and coroot lattices, respectively, and put $f=[\Pcoweight:\Qcoroot]$.
Thus, roots are integer valued on $\Pcoweight$, and the finite Weyl group $W$
acts trivially on $\Pcoweight/\Qcoroot$.  Set
$\rho^\vee=\sum_i\omega_i^\vee$, so that
$\alpha_i(\rho^\vee)=1$ and $\theta(\rho^\vee)=\hcox-1$.

\begin{lemma}[{\cite[Theorem 2.2]{AthanasiadisCatalan}}]\label{lem:band-modular-count}
Let $\mathcal B$ be a finite arrangement defined by integral affine
equations on a lattice $\mathsf L$ of rank $r$.  There is an integer $D_0>0$
such that, for every sufficiently large positive integer $q$ satisfying
$\gcd(q,D_0)=1$, its complement in $\mathsf L/q\mathsf L$ has
$\chi(\mathcal B,q)$ points.  The integer $D_0$ can be chosen using
only the coefficient matrices of subsets of the equations.
For type $A$ root normals on the full coweight lattice, one may take $D_0=1$. 
\end{lemma} 



Let $N(q)$ denote the number of points of $\Pcoweight/q\Pcoweight$ avoiding all the
equations in~\eqref{eq:band-definition} modulo $q$.

\begin{lemma}\label{lem:band-shrinking}
For every integer $q>M\hcox$, with $t=q-M\hcox$,
\begin{equation}\label{eq:band-shrunken-count}
 fN(q)=\sum_{w\in W_0}|\Pcoweight\cap(p+w(t\overline{\mathfrak a}))|.
\end{equation}
\end{lemma}

\begin{proof}
The quotient $\Pcoweight/q\Qcoroot$ is an $f$-fold cover of $\Pcoweight/q\Pcoweight$.
For $x\in \Pcoweight$, its class avoids the band equations precisely when
\begin{equation}\label{eq:band-modular-distance}
 |\alpha(x-p)-kq|\geq M
       \qquad(\forall\alpha\in\mathring{\Phi},\ k\in\mathbb Z).
\end{equation}
Indeed, $\alpha : \Pcoweight \to \mathbb{Z}$ is surjective and $\alpha (x-q\varpi) - \alpha (p) =\alpha (x-p) - q\alpha (\varphi)$, and the desired identity follows immediately.

Put $y=x-p$.  
Since $M>0$, every $y$ with $x$ satisfying
\eqref{eq:band-modular-distance} lies in a unique open alcove of the affine arrangement $\{ H_{\alpha ,kq}: k\in \mathbb{Z}\}$.  
Its alcoves modulo $q\Qcoroot$ have representatives $w(q\mathfrak a)$, for $w\in W_0$.
Inside $q\mathfrak a$, condition~\eqref{eq:band-modular-distance} is equivalent
to
\begin{equation}\label{eq:band-simple-inequalities}
 \alpha_i(y)\geq M\quad(1\leq i\leq r),
 \qquad q-\theta(y)\geq M.
\end{equation}
Necessity follows by taking the simple roots and the highest root.
For sufficiency, each positive root $\alpha$ is a nonnegative
integral combination of simple roots, and so is $\theta-\alpha$.
Thus $M\leq\alpha(y)\leq q-M$, which gives
\eqref{eq:band-modular-distance} for every root and every $k$.

The polytope in~\eqref{eq:band-simple-inequalities} is exactly
$M\rho^\vee+t\overline{\mathfrak a}$.  The allowed portions of distinct
$q$-alcoves are disjoint.  Consequently,
\begin{align*}
 fN(q)
 &=\sum_{w\in W_0}
   |(\Pcoweight-p)\cap w(M\rho^\vee+t\overline{\mathfrak a})|\\
 &=\sum_{w\in W_0}|\Pcoweight\cap(p+w(t\overline{\mathfrak a}))|,
\end{align*}
where the second equality uses $Mw(\rho^\vee)\in \Pcoweight$.
\end{proof}

\begin{lemma}\label{lem:band-stabilizers}
For every positive integer $t$ and every $x\in \Pcoweight$, there is an
isomorphism
\begin{equation}\label{eq:band-stabilizer-isomorphism}
 (t\Qcoroot\rtimes W_0)_{x-p}
       \cong\operatorname{Stab}_{W_F}(x\bmod t\Qcoroot).
\end{equation}
In particular, for $q>M\hcox$ and $t=q-M\hcox$,
\begin{equation}\label{eq:band-burnside-count}
 fN(q)=\sum_{g\in W_F}|(\Pcoweight/t\Qcoroot)^g|
\end{equation}
\end{lemma}

\begin{proof}
Suppose $(t\lambda,w) :z\mapsto wz+t\lambda$, with $\lambda\in \Qcoroot$, fixes
$y=x-p$.  Then
\[
 p-wp=x-wx-t\lambda\in \Qcoroot.
\]
It follows that $g(z)=wz+p-wp$ lies in the affine Weyl group, fixes $p$, and satisfies
$g(x)-x=-t\lambda$.  
Conversely, for every such $g = (\lambda',w)\in W_F$ fixing $x$ modulo $t\Qcoroot$, $g(x) -x=\lambda' + wx -x = -t\lambda \in t\Qcoroot$.
Then $(t\lambda ,w)$ fixes $y$.
These homomorphisms are inverse to each other and \eqref{eq:band-stabilizer-isomorphism} follows.

The number of closed $t$-alcoves containing a point $y$ equals
$|(t\Qcoroot\rtimes W_0)_y|$: the reflections in walls through $y$ generate
its finite stabilizer and act simply transitively on the incident alcoves. 
This shows that on the right-hand side of
\eqref{eq:band-shrunken-count}, elements in $x+ \Qcoroot$ appear $|(t\Qcoroot\rtimes W_0)_{x-p}|$-times.   
Therefore, the right-hand side of
\eqref{eq:band-shrunken-count} is equal to
\[
 \sum_{x+t\Qcoroot\in \Pcoweight/t\Qcoroot}|(t\Qcoroot\rtimes W_0)_{x-p}|
 =\sum_{x+t\Qcoroot\in \Pcoweight/t\Qcoroot}|\operatorname{Stab}_{W_F}(x+t\Qcoroot)|=\sum_{g\in W_F}|(\Pcoweight/t\Qcoroot)^g|.
\]
We are done.
\end{proof}

The arrangement $\mathcal{B}_M(p)$ is independent on the choice of $p$ in $F$.  
Indeed, for each $\alpha\in\mathring{\Phi}^+$, the allowed $j$'s such that $H_{\alpha ,j}\in\mathcal{B}_M(p)$ are
\begin{equation}\label{eq:interval_def0}
 \begin{cases}
  -M+1,\ldots,M-1,&\alpha(p)=0,\\
  -M+1,\ldots,M,&0<\alpha(p)<1,\\
  -M+2,\ldots,M,&\alpha(p)=1,
 \end{cases}
\end{equation}
and these ranges are constant on $F$.
We may thus assume that $p$ is a rational point of $F$ and choose an integer $N>0$ such that $Np\in \Pcoweight$.

\begin{lemma}\label{lem:band-fixed-points}
There is an integer $D_1>0$ such that, for every positive integer $t$ coprime to $D_1$ and every $g=(\alpha , w)\in W_F$
\begin{equation}\label{eq:band-fixed-point-number}
 (\Pcoweight/t\Qcoroot)^g
       =f\,t^{\dim\ker(w-I)}.
\end{equation}
\end{lemma}

\begin{proof}
Let $D_1$ be the least common multiple of $fN$ and all nonzero Smith factors of the matrices $w-I(w\in W_0)$ acting on $\Pcoweight$.
Since $\gcd(t,f)=1$, the natural map
\[
 \Pcoweight/t\Qcoroot\longrightarrow(\Pcoweight/\Qcoroot)\times(\Pcoweight/t\Pcoweight),
 \qquad x\longmapsto(x\bmod \Qcoroot,x\bmod t\Pcoweight),
\]
is an isomorphism.  

The element $g$ acts trivially on $\Pcoweight/\Qcoroot$.
Since $\alpha =p-w\alpha$ and $N$ is invertible in $\mathbb{Z}/t\mathbb{Z}$, $\alpha\in y'-w(y') + t\Pcoweight$ for some $y'\in \Pcoweight$. 
Therefore, $x+t\Pcoweight\in (\Pcoweight /t\Pcoweight)^g$ if and only if $w(x-y')-(x-y')\in t\Pcoweight$.
As $t$ being coprime to all nonzero Smith factors, we have $(w-I)(x-y') = 0$.
The kernel of $w-I|_{\Pcoweight}$ is pure $\mathbb{Z}$-submodule of $\Pcoweight$, there are $t^{\dim \ker (w-I)}$ points in $\Pcoweight/t\Pcoweight$ fixed by $g$. 
Combining with $f = |\Pcoweight /\Qcoroot|$, the identity \eqref{eq:band-fixed-point-number} follows.
\end{proof}

\begin{lemma}\label{lem:band-compatible-moduli}
There are infinitely many positive integers $q$ for which
\[
 \gcd(q,D_0)=1,
 \qquad \gcd(q-M\hcox,D_1)=1,
\]
where $D_0$ is chosen as in \cref{lem:band-modular-count} for
$\mathcal B_M(p)$ and $D_1$ as in
\cref{lem:band-fixed-points} such that $D_0 =1$ if $\mathring{\Phi}$ is of type $A_{2k}$.
\end{lemma}

\begin{proof}
For each odd prime $\wp$ dividing $D_0D_1$, choose a residue for $q$ avoiding $0$ if $\wp\mid D_0$ and avoiding $M\hcox$ if
$\wp\mid D_1$.  At most two residues are excluded, so this is
possible.  The same choice is possible at $\wp=2$ if $M\hcox$ is even.
If $M\hcox$ is odd, then $\hcox$ is odd.  The irreducible crystallographic
root systems with odd Coxeter number are precisely $A_{2k}$.
In that case \cref{lem:band-modular-count} allows $D_0=1$,
so there is at most one excluded residue modulo $2$.
The Chinese remainder theorem now produces a suitable arithmetic
progression, containing infinitely many positive integers.
\end{proof}

\begin{proof}[Proof of \cref{thm:affine-band-factorization}]
Choose sufficiently large $q$ from
\cref{lem:band-compatible-moduli} and put $t=q-M\hcox>0$.
\cref{lem:band-modular-count}, \cref{lem:band-stabilizers},
and \cref{lem:band-fixed-points} give
\[
 \chi(\mathcal B_M(p),q)
     =\sum_{g=(\alpha ,w)\in W_F}t^{\dim\ker(w-I)}.
\]
The Shephard--Todd--Solomon fixed-space identity for a finite real
reflection group is
\[
 \sum_{g\in W_F}t^{\dim\ker(w_g-I)}
      =t^{r-r_F}\prod_{i=1}^{r_F}(t+e_i).
\]
Here the power $t^{r-r_F}$ records the common fixed subspace.
This is the fixed-space identity of Solomon \cite[equation~(1)]{Solomon}.
Substituting $t=q-M\hcox$ proves
\eqref{eq:band-characteristic-polynomial} for infinitely many $q$,
and hence proves it as a polynomial identity.
\end{proof}

\begin{corollary}\label{cor:band-scopes-count}
The number of chambers of $\mathcal B_M(p)$ is
\begin{equation}\label{eq:band-region-count}
 (M\hcox+1)^{r-r_F}\prod_{i=1}^{r_F}(M\hcox+1-e_i).
\end{equation}
For an untwisted affine highest-weight module $V(\Lambda)$, the number of Scopes
classes in the orbit of $\Lambda-M\delta$ is
\[
 \frac{(M\hcox+1)^{r-r_F}\prod_{i=1}^{r_F}(M\hcox+1-e_i)}{|W_\Lambda|}.
\]
\end{corollary}

\begin{proof}
Apply the real-region formula
$|\operatorname{Chambers}(\mathcal A)|=(-1)^r\chi(\mathcal A,-1)$ to \cref{thm:affine-band-factorization}. For the number of Scopes classes,
use \cref{thm:band-decomposition} and the free stabilizer
action of \cref{thm:all-chambers}.
\end{proof}

The theorem proves the characteristic-polynomial factorization
uniformly.  No assertion of logarithmic freeness is needed in its
proof or in the chamber-count consequence.

\section{Factorization at arbitrary level}\label{sec:factorization}

\subsection{The depth shift in type C at every level}


Recall from \cref{sec:ArrangementAndNormalized}, let $a\in \mathbb{Z}^e$ s.t. $\ell \ge a_1\ge \ldots \ge a_e \ge 0$. 
For $\gamma\in \mathring{\Phi}$, set $d_\gamma = d_\gamma^+$ if $\gamma > 0$ and $d_\gamma =d^-_{-\gamma}$ if $\gamma <0$. 
Recall from \cref{prop:walls} that $b_\gamma (M) = M - d_\gamma$. 
Combining \cref{thm:band-decomposition} and \eqref{eq:interval_def0}, we have
\begin{equation}\label{eq:increment-bounds}
    M-1 \le b_{\gamma}(M+1) \le M+2, \quad b_{\gamma}(M+1) + b_{-\gamma} (M+1) \ge 2M.
\end{equation}

\begin{lemma}\label{lem:absent-root-closure}
Let $V$ be an integrable highest weight module of a complex symmetrizable Kac--Moody algebra $\mathfrak{g}$ and \(\eta\in P(V)\).  
Suppose that \(\alpha,\beta,\alpha+\beta\) are real roots.  If \(\eta+\alpha, \eta+\beta\not\in P(V)\), then
\(\eta+\alpha+\beta\not\in P(V)\).
\end{lemma}
\begin{proof}
Since $[\mathfrak{g}_\alpha ,\mathfrak{g}_\beta] =\mathfrak{g}_{\alpha + \beta}$, We have $\mathfrak{g}_{\alpha + \beta}\cdot V_{\eta} =0$ and $(\eta ,\rho) \ge 0$ for $\rho\in\{\alpha,\beta ,\alpha +\beta\}$.
The integrability of $V$ implies that $\mathfrak{g}_{-\alpha -\beta}\ni [e_{-\alpha},e_{-\beta}]: V_{\eta + \alpha + \beta} \to V_{\eta}$ is injective.
In particular, if $\eta + \alpha + \beta\in P(V)$, then either $e_{-\alpha}|_{V_{\eta +\alpha + \beta}}$ or $e_{-\beta}|_{V_{\eta +\alpha + \beta}}$ is non-zero, this contradicts to $\eta +\alpha,\eta +\beta \not\in P(V)$.
\end{proof}

\begin{lemma}\label{lem:strict-triangle}
Suppose that $\alpha,\beta,\alpha +\beta\in\mathring{\Phi}$ and $M \ge 1$, then
\begin{equation}\label{eq:strict-triangle}
 b_{\alpha + \beta}(M+1) < b_{\alpha}(M+1) + b_{\beta}(M+1).
\end{equation}
\end{lemma}
\begin{proof}
Let $\nu = \ell\Lambda_0 +a - \zeta_{\Lambda}(a)\delta$. 
Then, the real roots
\[\alpha' = \alpha + (1-d_\alpha)\delta \quad \text{and} \quad \beta' = \beta + (1-d_\alpha)\delta\]
satisfies $\nu + \alpha',\nu+\beta'\not\in P(\Lambda)$.
By \cref{lem:absent-root-closure}, $\nu +\alpha'+\beta'\not\in P(\Lambda)$, whence
\[
d_{\alpha} + d_{\beta} -2 <d_{\alpha + \beta}.
\]
Equivalently, $2M - b_{\alpha}(M+1) - b_\beta(M+1)\le M - b_{\alpha + \beta}(M+1)$. The desired inequality follows since $M > 0$.
\end{proof}

Recall that $W_0$ is the group of signed permutation acting on $\mathring{\mathfrak{h}}^*_{\mathbb{R}}$.
For $\sigma \in W_0$, define
\begin{equation}\label{eq:gap-statistic}
    A_M(\sigma)=\sum_{i=1}^{e-1}b_{\sigma (\alpha_i)}(M+1)
       +\left\lceil\frac{b_{-\sigma (\theta)}(M+1)}{2}\right\rceil
       +\left\lfloor\frac{b_{\sigma (\alpha_e)}(M+1)}{2}\right\rfloor.
\end{equation}
 

\begin{lemma}
\label{lem:gap-count}
For every \(M\ge1\), one has the polynomial identity
\begin{equation}\label{eq:gap-count}
 \chi_{a,M}(q)=\sum_{\sigma \in W_0}
  \binom{(q-1)/2-A_M(\sigma)+e}{e},
\end{equation}
where 
the binomial coefficient denotes its degree-\(e\) polynomial.
\end{lemma}
\begin{proof}

We apply \cref{prop:finite-field} to the normalized arrangement of $\mathcal{A}_{a,M}$ to compute $\chi_{a,M}(q)$ for primes $q =2Q+1\gg 0$. 
Every point $p\in \mathbb{F}_q^e$ has a unique representative $p' =(p'_1,\ldots ,p'_e)\in \mathbb{Z}^e$ such that $|p'_i|\le Q$ and $p$ avoids all walls if and only if so does $p'$. 
In this case, $p'$ is uniquely determined by $0\le k_1 \le k_2 \le \ldots \le k_e \le Q$ and $\sigma \in W_0$ such that $(k_1,\ldots ,k_e) = \sigma^{-1} (p')$ and 
\begin{equation}\label{eq:gap_criteria}
    \begin{aligned}
        & 2k_1 \ge b_{-\sigma (\theta)}(M+1), \\
        & k_{i+1}-k_{i} \ge b_{\sigma(\alpha_i)}(M+1) \quad (1\le i < e), \\
        & q-2k_e \ge b_{\sigma (\alpha_e)}(M+1).
    \end{aligned}
\end{equation}
Indeed, $(k_1,\ldots , k_e)$ is obtained by arranging $(|p'_1|,\ldots , |p'_e|)$ in non-decreasing order and $\sigma$ is unique by \eqref{eq:increment-bounds}.
Conversely, given $0\le k_1\le \ldots \le k_e \le Q$ and $\sigma$ such that \eqref{eq:gap_criteria} holds, then $p'\in \sigma (k_1,\ldots ,k_e)$ avoids the walls in $\mathbb{F}_q^e$ by repeatedly using \eqref{eq:strict-triangle}.

Fixing $w\in W_0$, we compute the number of sequences $0\le k_1 \le \ldots \le k_e \le Q$ that satisfy \eqref{eq:gap_criteria}. 
Set 
\[
    k'_i = k_i - \left\lceil\frac{b_{-\sigma (\theta)}(M+1)}{2}\right\rceil - \sum_{j=1}^{i-1} b_{\sigma (\alpha_i)}(M+1),
\]
then \eqref{eq:gap_criteria} is equivalent to
\[
0\le k'_1 \le k'_2 \le \ldots \le k'_e \le Q - A_M(\sigma).
\]
So, the number of such sequences is $\binom{(q-1)/2-A_M(\sigma)+e}{e}$. Summing over all $w\in W_0$ yields the identity \eqref{eq:gap-count}.
\end{proof}

\begin{theorem}\label{thm:depth-shift}
For every admissible $a$, every rank \(e\ge2\),
and every depth \(M\ge1\),
\begin{equation}\label{eq:depth-shift-proved}
 \chi_{a,M}(q)=\chi_{a,1}\bigl(q-2e(M-1)\bigr).
\end{equation}
In particular, integral factorization at all positive depths is
equivalent to integral factorization at depth one.
\end{theorem}
\begin{proof}

If $M$ is odd, $A_M(\sigma) =A_1(\sigma)+e(M-1)$. 
If $M$ is even, 
\begin{align*}
    A_M(\sigma) & =\sum_{i=1}^{e-1}b_{\sigma (\alpha_i)}(2)+e(M-1) + \left\lceil \frac{b_{-\sigma (\theta)}(2)+1}{2} \right\rceil + \left\lfloor \frac{b_{\sigma (\alpha_e)}(2)-1}{2} \right\rfloor \\
    & = \sum_{i=1}^{e-1}b_{\sigma (\alpha_i)}(2)+e(M-1) + \left\lfloor \frac{b_{-\sigma (\theta)}(2)}{2} \right\rfloor + \left\lceil \frac{b_{\sigma (\alpha_e)}(2)}{2} \right\rceil \\
    & = A_1(\sigma^*) + e(M-1).
\end{align*}
where $\sigma^*\in W_0$ is such that $\sigma^* (i) = -\sigma (e+1-i)$.
The identity then follows from \eqref{eq:gap-count}.
\end{proof}

\subsection{The level-one specialization}
The general zero-content theorem below includes level one. We record
its specialization for the normal forms of \cref{sec:counting}.
Write $\operatorname{Exp}(C_n)=\{1,3,\ldots,2n-1\}$, empty for $n=0$.

\begin{corollary}\label{thm:level-one-factorization}
Fix \(0\le k,t\le e\), with \(t\equiv k\pmod2\), and
\(M\ge1\). For the orbit with \(h(p)=t\) in \(V(\Lambda_k)\),
the roots of \(\chi_M(q)\), counted with multiplicity, are
\[
\begin{cases}
 \{2eM-d:d\in\operatorname{Exp}(C_t)\}
 \sqcup\{2eM-2+d:d\in\operatorname{Exp}(C_{e-t})\},&t<k,\\
 \{2eM-d:d\in\operatorname{Exp}(C_t)\sqcup \operatorname{Exp}(C_{e-t})\},&t=k,\\
 \{2eM-2+d:d\in\operatorname{Exp}(C_t)\}
 \sqcup\{2eM-d:d\in\operatorname{Exp}(C_{e-t})\},&t>k.
\end{cases}
\]
\end{corollary}
\begin{proof}
This corollary follows from \cref{thm:zero-content-factor}, whose proof is independent of this corollary.  
Here we give a direct proof.
Applying the diagram involution $i\mapsto e-i$, the result for $t>k$ can be derived from that of $t<k$, so we may assume that $t\le k$. Using the table in the proof of \cref{prop:vertex-bands}, we can show that
\begin{align*}
    & \left\lceil \frac{b_{-\sigma (\theta)}(2)}{2} \right\rceil = \begin{cases}
        0 & \text{if } -\sigma (1) \in [1,t], \\
        1 & \text{otherwise,}
    \end{cases} \quad \left\lfloor \frac{b_{-\sigma (\alpha_e)}(2)}{2} \right\rfloor = \begin{cases}
        1 & \text{if } -\sigma (e) \in [1,t], \\
        0 & \text{otherwise,}
    \end{cases} \\
    & b_{\sigma (\alpha_i)}(2) = \begin{cases}
        0 & \text{if } -\sigma (i+1), \sigma (i)  \in [1,t] ,\\
        2 & \text{if } -\sigma (i) \in [1,t]\text{ and } -\sigma (i+1)\not\in [1,t], \\ 
        2 & \text{if } |\sigma (i)| >t \text{ and } \sigma (i+1)\in [1,t] \\
        2 & \text{if } \min (|\sigma (i)|,|\sigma (i+1)|)>t \text{ and } t < k, \\
        1 & \text{otherwise},
    \end{cases} 
\end{align*}
So,
\[
A_1(\sigma) = \begin{cases}
    e + s + \delta_{t<k}(e - t - s-1)  & \text{if } |\sigma (1)|,|\sigma (e)| > t, \\
    e + s-1 + \delta_{t < k}(e - t -s)  & \text{if } |\sigma (1)|\le t, |\sigma (e)| > t, \\
    e + s + \delta_{t < k}(e - t -s) & \text{if } |\sigma (1)|>t, |\sigma (e)| \le t, \\
    e + s-1 + \delta_{t < k}(e - t -s+1) & \text{if } |\sigma (1)|,|\sigma (e)| \le t,
\end{cases}
\]
where $s$ is the number of intervals in $T(\sigma) = \{ 1\le i \le t : |\sigma (i)| \le t \}$. 
In particular, $A_1(\sigma)$ is a constant on any right $W_t\times W_{e-t}$-coset.
If $k=t$, then
\begin{align*}
    \frac{\chi_{a,1}(q)}{2^et!(e-t)!}=&  \sum_{s=0}^{t} \Bigg(\binom{t-1}{s-1}\left(\binom{e-t-1}{s} + \binom{e-t-1}{s-1}\right)\binom{Q - s}{e} \\
    & \quad + \binom{t-1}{s-1}\left(\binom{e-t-1}{s-1} + \binom{e-t-1}{s-2}\right)\binom{Q -s+1}{e}\Bigg) \\
    =& \sum_{s=0}^{t} \binom{t-1}{s-1}\Bigg(\binom{e-t}{s}\binom{Q -s}{e}  + \binom{e-t}{s-1}\binom{Q -s+1}{e}\Bigg) \\
    = & \sum_{s = -1}^{t}\left(\binom{t-1}{s-1}+\binom{t-1}{s}\right)\binom{e-t}{s}\binom{Q -s}{e} \\
    = & \sum_{s=0}^t \binom{t}{s}\binom{e-t}{s}\binom{Q-s}{e} \\
    = & \sum_{s =0}^t\sum_{r=0}^{t-s}  \binom{t}{s}\binom{e-t}{s}\binom{Q-t}{e-r}\binom{t-s}{r} \\
    = & \sum_{r = 0}^t \binom{Q-t}{e-r}\sum_{s=0}^{t-r}\binom{t}{s}\binom{e-t}{s}\binom{t-s}{r} \\
    = & \sum_{r = 0}^t \binom{Q-t}{e-r}\binom{t}{r}\sum_{s=0}^{t-r}\binom{t-r}{s}\binom{e-t}{e-t-s} \\
    = & \sum_{r = 0}^t \binom{Q-t}{e-r}\binom{t}{r}\binom{e-r}{e-t} \\
    = & \binom{Q-t}{e-t}\sum_{r=0}^t\binom{Q-e}{t-r}\binom{t}{r} \\
    = & \binom{Q-t}{e-t} \binom{Q-e+t}{t}.
\end{align*}
where $Q = (q-1)/2$. So $\chi_{a,1}(q) = \prod_{s=1}^{e-t}(q - (2e-2s+1))\prod_{r=1}^t(q-(2e-2r+1))$.
The computation for $t < k$ is simpler and we leave it to the readers.
\end{proof}

\subsection{Upward roots when all content coefficients are positive}

\begin{lemma}\label{lem:short-upward-maximal}
Let $\Lambda\in P^+$ be a dominant integral weight of an untwisted affine
root system, with $\nu\in P^+(\Lambda)$ and $\nu+\delta\notin P(\Lambda)$. 
Write $\alpha=\Lambda-\nu\in Q_+$. 
If $\beta$ is a positive short real root, then
\[
 \nu+\beta\in P(\Lambda)
 \quad\Longleftrightarrow\quad \beta\leq \alpha.
\]
\end{lemma}
\begin{proof}
Necessity follows from $P(\Lambda)\subset \Lambda -Q_+$. 
For sufficiency, assume $\gamma=\alpha-\beta\in Q_+$. 
Since $\nu$ is maximal and dominant, we have $\alpha\not\ge \delta$ by \cite[Proposition 11.3(a)]{K90}.
The same argument yields that $\nu + \beta\in P(\Lambda)$ if $\nu + \beta\in P^+$.

Suppose $\nu+\beta \not\in P^+$, and choose a simple root $\alpha_i$ with $\langle\nu+\beta,\alpha_i^\vee\rangle<0$. 
Then $-2 \le \langle \beta ,\alpha_i^\vee \rangle <0$ since $\beta$ is short. 
The equality holds if and only if $\beta + \alpha_i = n\delta$ for some $n>0$.
In this case, $\langle \gamma ,\alpha_i^\vee\rangle = \langle \Lambda ,\alpha_i^\vee \rangle - \langle \nu + \beta , \alpha^\vee_i\rangle >0$. 
But this implies that $\gamma \ge \alpha_i$ and thus $\alpha \ge \beta +\alpha_i = n\delta$, a contradiction.
So, we must have $\langle \beta ,\alpha_i^\vee \rangle =-1$ and $\langle \nu ,\alpha^\vee_i \rangle =0$.
Now, $s_i (\nu + \beta) = \nu + \beta + \alpha_i = \nu + \beta'$, $\beta'$ is a positive short real root and $\operatorname{ht}(\beta') > \operatorname{ht}(\gamma)$. 

Repeat this operation whenever the resulting weight is not dominant.
The height of the root increases strictly and stays bounded by the height
of $\alpha$, so the procedure terminates at a dominant weight $\nu+\beta'$ with $\beta'\leq \alpha$ which has been proved.
\end{proof}

\begin{proposition}
\label{prop:positive-content-upward}
Work in type $C_e^{(1)}$, with
\[
 \alpha_0=\delta-2\varepsilon_1,\qquad
 \alpha_i=\varepsilon_i-\varepsilon_{i+1}\ (1\leq i<r),\qquad
 \alpha_e=2\varepsilon_e.
\]
Let $\nu\in P^+(\Lambda)$ be a maximal weight, and write $\Lambda-\nu=\sum_{i=0}^e c_i\alpha_i$. Suppose every $c_i$ is positive.
Then the set $\{i:1\leq i<e,\ c_i=1\}$ is nonempty; write its first and
last indices as $L$ and $R$. For every positive real root $\beta$,
\[
 \nu+\beta\in P(\Lambda)
 \quad\Longleftrightarrow\quad \beta\leq\Lambda-\nu.
\]
More explicitly, the upward roots are precisely
\[
\begin{array}{c|c}
\text{root}&\text{condition}\\\hline
\varepsilon_i-\varepsilon_j&i<j\\
\varepsilon_i+\varepsilon_j&i<j,\ j>R\\
\delta-\varepsilon_i-\varepsilon_j&i<j,\ i\leq L\\
\delta-\varepsilon_i+\varepsilon_j&i<j,\ i\leq L,\ j>R\\
2\varepsilon_i&i>R\\
\delta-2\varepsilon_i&i\leq L.
\end{array}
\]
\end{proposition}
\begin{proof}
We've seen in the proof of \cref{lem:short-upward-maximal} that $\alpha =\Lambda -\nu\not\ge \delta = \alpha_0 + 2\alpha_1 + \ldots + 2\alpha_{e-1}+\alpha_e$. 
The first statement follows.

Every positive real root below $\alpha$ is either $\gamma$ or $\delta-\gamma$
for $\gamma\in\mathring{\Phi}^+$. 
The simple-root expansions
\[
\begin{aligned}
\varepsilon_i-\varepsilon_j
 &=\alpha_i+\cdots+\alpha_{j-1},\\
\varepsilon_i+\varepsilon_j
 &=\alpha_i+\cdots+\alpha_{j-1}
      +2(\alpha_j+\cdots+\alpha_{e-1})+\alpha_e,\\
\delta-\varepsilon_i-\varepsilon_j
 &=\alpha_0+2(\alpha_1+\cdots+\alpha_{i-1})
      +\alpha_i+\cdots+\alpha_{j-1},\\
\delta-\varepsilon_i+\varepsilon_j
 &=\alpha_0+2(\alpha_1+\cdots+\alpha_{i-1})
      +\alpha_i+\cdots+\alpha_{j-1}
      +2(\alpha_j+\cdots+\alpha_{e-1})+\alpha_e
\end{aligned}
\]
give exactly the four short-root conditions in the table. Lemma~
\ref{lem:short-upward-maximal} proves support for every such candidate.
The two long-root expansions are
\[
 2\varepsilon_i=2(\alpha_i+\cdots+\alpha_{e-1})+\alpha_e,
 \qquad
 \delta-2\varepsilon_i=\alpha_0+2(\alpha_1+\cdots+\alpha_{i-1}),
\]
so their necessary conditions are respectively $i>R$ and $i\leq L$.

We prove sufficiency for $2\varepsilon_i$ with $i>R$. 
Write $\nu = \ell \Lambda_0 + a -c_0\delta$ with $\ell \ge a_1 \ge \ldots \ge a_e\ge 0$.
We have
\[
 \langle\nu,\alpha_R^\vee\rangle
 =\langle\Lambda,\alpha_R^\vee\rangle - a_{R,R-1}c_{R-1} -2 - a_{R,R+1}c_{R+1} \ge 1,
\]
where we use $- a_{R,R+1}c_{R+1} \ge 2$ since $-a_{R,R+1} = 2$ if $R=e-1$ and $c_{R+1}=2$ if $R< e-1$.

Consider $\nu+2\varepsilon_i$. If $a_{i-1}-a_i\geq2$, it is already
dominant. If $a_{i-1}-a_i=1$, reflection in $\alpha_{i-1}$ changes it to
\[
 \nu+\varepsilon_{i-1}+\varepsilon_i.
\]
This short root is below $\alpha$ because $i>R$, so it belongs to $P(\Lambda)$ by \cref{lem:short-upward-maximal}, and so does $\nu + 2\varepsilon_i$. 
In particular, $\nu + 2\varepsilon_R \in P(\Lambda)$.
If $a_{i-1}=a_i$, the reflection fixes $\nu$ and changes the long root to $2\varepsilon_{i-1}$. 
As long as $i-1>R$, this remains below $\alpha$, and the argument repeats. 
It cannot pass from index $R+1$ to $R$ through an equality, because $a_R-a_{R+1}\geq1$. 
Thus the process terminates in one of the preceding cases. 

Finally, applying the Dynkin diagram involution $j \mapsto e-j$, we obtain the result for $\delta - 2\varepsilon_i$.
\end{proof}

\begin{corollary}
\label{cor:positive-content-thresholds}
Under the hypotheses of \cref{prop:positive-content-upward},
put \(I_-=\{1,\ldots,L\}\), \(I_0=\{L+1,\ldots,R\}\), and
\(I_+=\{R+1,\ldots,r\}\). Then
\[
\begin{aligned}
 d_{\eps_i-\eps_j}^+&=0,&
 d_{\eps_i-\eps_j}^-&=-\delta_{i\in I_-,\ j\in I_+}
 &&(i<j),\\
 d_{\eps_i+\eps_j}^+&=\delta_{j\in I_-\cup I_0},&
 d_{\eps_i+\eps_j}^-&=-\delta_{i\in I_-}
 &&(i\le j),
\end{aligned}
\]
where \(i=j\) denotes the long root \(2\eps_i\) and $\delta_{P}$ is the Kronecker delta.
\end{corollary}
\begin{proof}
Recall that $d^{\pm}_\gamma$ is the smallest integer such that
\[
\nu \pm \gamma -d^{\pm}_{\gamma}\delta \in P(\Lambda).
\]
We have proved that $a_R > a_{R+1}$. 
We can show that $a_L > a_{L+1}$ in the same way.

Fix $\gamma \in \mathring{\Phi}^+$. 
By \cref{prop:wall-growth}, $d^+_\gamma \ge 0$ and $d^-_\gamma \ge -1$.
If $(\nu,\gamma)> 0$, then $\nu - \gamma \in P(\Lambda)$ and hence $d^-_\gamma \le 0$.

If $(\nu, \gamma)=0$, then $a_i = a_j$ and $\gamma = \varepsilon_i -\varepsilon_j(i<j)$ or $a_i = a_j =0$ and $\gamma =\varepsilon_i+\varepsilon_j(i\le j)$.
In the latter case, $i,j>R$.
Then $\nu -\gamma \in P(\Lambda)$ since $\nu +\gamma =s_\gamma (\nu -\gamma) \in P(\Lambda)$ by \cref{prop:positive-content-upward}. 
This proves that $d^-_\gamma \le 0$.

Likewise, if $(\delta -\gamma,\nu) >0$ gives $d_\gamma^+ \le 1$. 
If $(\delta -\gamma,\nu) =0$, then $\gamma = \varepsilon_i +\varepsilon_j$ with $i \le j$ and $a_i = a_j =\ell$.
In particular, $j\le L$.
By \cref{prop:positive-content-upward}, $\nu + \delta -\gamma \in P(\Lambda)$ and so is $\nu + \gamma - \delta =s_{\delta -\gamma}(\nu + \delta -\gamma)$. Thus, $d^+_\gamma \le 1$.




Consequently \(d_\gamma^+\in\{0,1\}\), with zero exactly when
\(\gamma\) is an upward root, and
\(d_\gamma^-\in\{-1,0\}\), with minus one exactly when
\(\delta-\gamma\) is an upward root. The asserted formulas
now follow from the six root families in the proposition.
\end{proof}

\begin{theorem}\label{thm:positive-content-factor}
Keep the assumption and notations in \cref{prop:positive-content-upward}.
Put \(r_A=L\), \(r_B=R-L\), \(r_C=e-R\), and
\(d=r_A+r_C\). Then, for every integer \(M\ge0\),
\begin{equation}\label{eq:positive-content-factor}
 \chi_M(q)=
 \prod_{u\in\operatorname{Exp}(C_{r_A})\sqcup\operatorname{Exp}(C_{r_C})}
       (q-2eM-2d+u)
 \prod_{j=0}^{r_B-1}(q-2eM-2d-j).
\end{equation}
In particular the factors depend only on \(L\) and \(R\).
\end{theorem}
\begin{proof}
Write \(I_-=\{1,\ldots,L\}\), \(I_0=\{L+1,\ldots,R\}\),
and \(I_+=\{R+1,\ldots,e\}\). 
By \cref{cor:positive-content-thresholds}, \(b_\gamma (1) \ge 0\), and $b_{\gamma}(1)=0$ if and only if $\gamma = \varepsilon_i + \varepsilon_j$ with $i\le j\le R$. 
Moreover, $b_{\gamma}(1) + b_{-\gamma} (1) \ge 1$ and $b_{\alpha + \beta}(1) \le b_\alpha (1) + b_\beta (1)$.

The proof of \cref{lem:gap-count} therefore works for \(M=0\).
Indeed, we only need the weak strictness of \eqref{eq:strict-triangle} in the proof.
Hence the same signed-permutation sum counts \(\chi_0\). 
Then \cref{thm:depth-shift} holds for $M=0$, that is,
\begin{equation}\label{eq:positive-zero-shift}
 \chi_M(q)=\chi_0(q-2eM)\qquad(M\ge0).
\end{equation}

For $M =0$, the arrangement on the \(I_-\cup I_+\) coordinates is exactly $B_1(p)$ in type $C_d^{(1)}$ $a=u=1/2\sum_{i=1}^{r_A}\varepsilon_i$.
\cref{thm:affine-band-factorization} gives its polynomial
\[
 \chi_{\rm core}(q)=
 \prod_{u\in\operatorname{Exp}(C_{r_A})\sqcup\operatorname{Exp}(C_{r_C})}
 (q-2d+u).
\]
On $\mathbb{F}_q^d$ with $q\gg 0$, if $x\in \mathbb{F}_q^d$ avoids all the wall, then
\[
 |\{x_i,1-x_i:i\in I_-\}\ \cup\ \{x_j,-x_j:j\in I_+\}| = 2d.
\]
Conversely, a point $x$ satisfying this condition must avoid all the walls.
By \cref{cor:positive-content-thresholds}, each \(I_0\) coordinate must avoid these \(2d\) values, and the
\(I_0\) coordinates must be pairwise distinct. There are no other
conditions on them. Thus
\[
 \chi_0(q)=\chi_{\rm core}(q)\prod_{j=0}^{r_B-1}(q-2d-j),
\]
first by finite-field counting and then as a polynomial identity.
Together with \eqref{eq:positive-zero-shift}, this proves the formula.
\end{proof}

\subsection{Upward roots in a finite Levi subsystem}

\begin{proposition}
\label{prop:zero-content-upward}
Let $\Lambda$ and $\nu$ be dominant integral weights of type $C_e^{(1)}$, with 
\[
 \Lambda-\nu=\sum_{i=0}^e c_i\alpha_i\in Q_+.
\]
Suppose $\min_i \{c_i\} =0$, and set
$J=\{i:c_i>0\}$. 
Then $N_\nu \subseteq \Phi_J$, and $\mu \in \nu + N_\nu$ is a weight of $U_q(\mathfrak{g}_J)v_\Lambda$, with $v_\Lambda$ being a highest weight vector.
We can treat $N_\nu$ according to the connected components of $J$.
An interior component of type $A$ contributes every one of its positive
roots.

For an end component of type $C_n$, choose its standard simple roots
\[
 \eta_i=\epsilon_i-\epsilon_{i+1}\quad(1\leq i<n),
 \qquad \eta_n=2\epsilon_n,
\]
and let $b_1,\ldots,b_n>0$ be the content coefficients on that component.
Let $a_1\geq\cdots\geq a_n\geq0$ be the finite coordinates of the
restricted dominant weight. Put
\[
 k=\max(\{i:1\leq i<n,\ b_i=1\} \cup \{0\}).
\]
The upward short roots are exactly
\[
 \epsilon_i-\epsilon_j\quad(i<j),
 \qquad
 \epsilon_i+\epsilon_j\quad(i<j,\ j>k).
\]
The upward long roots are as follows:
\begin{enumerate}
\item If $k\neq 1$, $2\epsilon_i\in N_{\nu}$ if and only if $i>k$.
\item If $k=1$, $2\epsilon_i\in N_\nu$ if and only if $i>t =\max \{ 1\le j \le n : a_j = a_1 \}$.
\end{enumerate}
For $n=1$, the convention gives $k=0$ and the sole positive root is
upward. A left end component is read in these coordinates after applying
the affine Dynkin involution $\alpha_i\mapsto\alpha_{r-i}$.
\end{proposition}
\begin{proof}
If $\nu+\beta\in P(\Lambda)$, then
$\beta\leq\Lambda-\nu$. 
Thus $\beta\in \Phi_J$. 
Because $J\subsetneq \{0,1,\ldots,e\}$, the root system $\Phi_J$ is finite. 
In particular, $\beta$ lies in a single connected component of $\Phi_J$.
Since $V(\Lambda)$ is a quotient of $U_q^-(\mathfrak{g})$ as a vector space, the PBW basis theorem implies that $\mu \in \nu + N_\nu$ must be a weight of $U_q(\mathfrak{g}_J)v_{\Lambda}$.
For disconnected $J$, it is the tensor product over the connected components. 

Let $\Lambda',\nu'\in P^+_J$ with $\nu' \in P(\Lambda')$.
We claim that if a positive real root $\beta\in \Phi_J$ satisfies $\beta \le \Lambda' -\nu'$, then $\nu' + \beta \in P(\Lambda')$.
The proof is analogue to \cref{lem:short-upward-maximal}.

Let $J_c$ be a component of $J$. 
If $0,e\not\in J_c$, then every positive root $\alpha$ is short and $\alpha \le \Lambda -\nu$, hence $\alpha \in N_\nu$ by the claim.

For an end component, the expansions
\[
\begin{aligned}
 \epsilon_i-\epsilon_j
   &=\eta_i+\cdots+\eta_{j-1},\\
 \epsilon_i+\epsilon_j
   &=\eta_i+\cdots+\eta_{j-1}
      +2(\eta_j+\cdots+\eta_{n-1})+\eta_n,\\
 2\epsilon_i
   &=2(\eta_i+\cdots+\eta_{n-1})+\eta_n
\end{aligned}
\]
show that the exact short-root conditions are those stated. They also
show that $i>k$ is necessary for an upward long root.

Suppose first that $k=0$. 
Then $b_i\ge 2, \forall i<n$, so $2\epsilon_i \le \Lambda' -\nu'$. 
Starting with $\nu'+2\epsilon_i$, compare $a_i+2$ to $a_{i-1}$ when $i>1$.
A gap $a_{i-1}-a_i\geq2$ makes the weight dominant. A gap equal to one
allows reflection in $\eta_{i-1}$ to give
$\nu'+\epsilon_{i-1}+\epsilon_i$, which is supported by the short-root
claim. A zero gap permits reflection fixing $\nu'$ and replacing
$2\epsilon_i$ by $2\epsilon_{i-1}$. All intermediate long roots remain $\le \Lambda'-\nu'$. 
This process either reaches one of the $\nu'+2\varepsilon_j$ with $a_{j-1} -a_j \ge 1$ case or reaches $\nu'+2\epsilon_1\le \Lambda'$, which is dominant. 
Thus every long positive root belong to $N_{\nu'}$. 
This also proves the $n=1$ case directly; equivalently, it is the elementary $\mathfrak{sl}_2$ string statement.

Next let $k\geq2$. Note that
\[
 a_k - a_{k+1} =\langle\nu',\eta_k^\vee\rangle= \langle \Lambda',\eta^\vee_k \rangle -a_{k,k-1}b_{k-1}-2-a_{k,k+1}b_k\ge -a_{k,k-1}b_{k-1} >0.
\]
The rest of proof is the same as the case $k=0$.

Finally, let $k=1$. 
If $i\le t$, then $\exists w\in W_{\nu'}$ such that 
$w(\nu' + 2\epsilon_i) = \nu' + 2\epsilon_1$, which is dominant.
But $\nu' + 2\epsilon_1\not\in P(\Lambda')$ as $b_1 =1$, so $i > t$ is necessary.
If $i>t$, the proof is the same as the case $k =0$.
This proves sufficiency and completes all cases.
\end{proof}

\begin{lemma}\label{lem:zero-content-boundary}
Keep the setting of \cref{prop:zero-content-upward}. 
Recall that $W_{\nu} = W_{J_\nu}$ is a standard parabolic subgroup of $W$.
Then each component of $J_\nu$ is contained in some component of $J$ or some component of $I\backslash J$. 
Moreover, if $i\in J_\nu \backslash J$, then $i\pm 1 \not\in J$.
\end{lemma}

\begin{proof}
If $c_i=0$, then
\[
 \langle\nu,\alpha_i^\vee\rangle
 =\langle\Lambda,\alpha_i^\vee\rangle
    -\sum_{|j-i|= 1} a_{ij}c_j \ge 0.
\]
The inequality is strict if $c_j > 0$ for some $j$ such that $|j-i| =1$.
In particular, $a_i\neq a_j$ whenever $i,j\in J$ belong to different components of $J$.
Recall that $i\in J_\nu$ if $a_i = a_{i+1}$, where $a_0 = \ell$ and $a_{e+1} =0$.
This completes the proof.

\end{proof}

\subsection{A weighted compression lemma}
In this subsection, $\nu = \ell \Lambda_0 + a -\zeta_\Lambda (a)\delta$. 
The Scopes arrangement under consideration is defined by $N_{\nu-\delta}$.
The strict triangle inequality \eqref{eq:strict-triangle} yields that
\begin{equation}\label{eq:weak-core-triangle}
 b_{\alpha + \beta}(1) \le b_\alpha (1) + b_\beta (1) \quad \forall \alpha,\beta,\alpha + \beta\in \mathring{\Phi}.
\end{equation}
In this subsection, we assume that $\min_i c_i =0$, then by \cref{lem:window}
\begin{equation}\label{eq:zero-window-core}
 b_\gamma (1) + b_{-\gamma}(1) = 1 + \kappa_{\gamma} \in \{0,1,2\}.
\end{equation}
Define an affine arrangement $\mathcal U$ and an affine reflection group $H$ by
\[
 \mathcal U=\{H_{\gamma ,d_\gamma}:\kappa_\gamma=1\},
 \qquad
 H=\langle s_{\gamma,b_{-\gamma}(1)}:\kappa_\gamma=-1\rangle,
\]
where $s_{\gamma,b}$ denotes reflection in $H_{\gamma ,b}$ and one uses one representative of each root direction. These definitions do not depend
on its sign.

\begin{theorem}\label{thm:weighted-core}
$H$ is a finite group and preserves $\mathcal U$. 
For all sufficiently large odd integers $q$, set $t=q-2e$. Then the number of points modulo $q$ outside the Scopes arrangement $N_{\nu - \delta}$ equals
\begin{equation}\label{eq:weighted-core-count}
 \sum_{X\in(\mathbb Z/t\mathbb Z)^e\setminus\mathcal U}|H_X|
 =\sum_{g\in H}|((\mathbb Z/t\mathbb Z)^e       \setminus\mathcal U)^g|.
\end{equation}
Consequently, its characteristic polynomial is obtained from the polynomial on the right by substituting $t=q-2e$.
\end{theorem}
\begin{proof}

For $s\in \{\pm1 ,\ldots ,\pm e\}$, define $h_s = (b_{\varepsilon_1 -\varepsilon_s}(1)-b_{\varepsilon_1+\varepsilon_s}(1))/2$ where $\varepsilon_{-s} = -\varepsilon_s$, and $b_0(1):=0$ by convention.
By \eqref{eq:weak-core-triangle} and \eqref{eq:zero-window-core}, we have
\begin{align*}
& b_{\varepsilon_1-\varepsilon_s}(1) - b_{\varepsilon_1- \varepsilon_m}(1) \ge - b_{\varepsilon_s - \varepsilon_m}(1),  \\
& -b_{\varepsilon_1+\varepsilon_s}(1) + b_{\varepsilon_1+ \varepsilon_m}(1) \ge -b_{\varepsilon_s- \varepsilon_m}(1), 
\end{align*}
hence, $F(s,m) =h_s - h_m +b_{\varepsilon_s -\varepsilon_m}(1)\ge 0$, equivalently, $(\gamma ,h) \ge -b_{\gamma}(1)$.
We have $F(s,m) + F(m,s) = \kappa_{\varepsilon_s -\varepsilon_m}+1$.
If $\kappa_{\gamma} =-1$, then $0=(\gamma ,h) + (-\gamma ,h) \ge -(1 + \kappa_{\gamma}) =0$, and thus $(\gamma,h)= -b_{\gamma}(1)$, that is, $h\in H_{\gamma,-b_{\gamma}(1)}$.
So, $H$ is a reflection subgroup of the finite reflection group $W_F$, where $F$ is the facet containing $h$.
Under the translated coordinate $y = x-h$, $H$ acts as signed permutations: let $\gamma = \varepsilon_s -\varepsilon_m$ with $\kappa_\gamma =-1$, for any $x\in \mathbb{R}^e$ with $y =x-h$, set $z =s_{\gamma,-b_{\gamma}(1)}(x)-h$, then $z_s = x_m-b_{\gamma}(1)-h_s = x_m -h_m = z_m$, similarly $z_m = y_s$ and $z_r = y_r$ for $|r|\not\in \{|s|,|m| \}$. 
The following weak triangle formulas
\[
-b_{\varepsilon_s-\varepsilon_m}(1) \le b_{\varepsilon_m - \varepsilon_r}(1) - b_{\varepsilon_s-\varepsilon_r}(1) \le b_{\varepsilon_m - \varepsilon_s}(1)
\]
thus the equalities hold and $F(m,r) = F(s,r)$. Similarly, $F(r,m) =F(r,s)$.

Suppose that $\alpha = \varepsilon_r - \varepsilon_{r'}$ with $\kappa_\alpha =1$, the $y$-coordinate of $H_{\alpha , d_\alpha} = H_{\alpha,b_{-\alpha}(1)-1}$ satisfies
\[
y_r - y_{r'} + 1-F(r',r) = 0.
\]
and it is fixed by $s_{\gamma,d_\gamma -1}$ unless $\{ r,r' \}\cap \{ s,m \} \neq \emptyset$. 
Without loss of generality, assume that $r'=s$, then $s_{\gamma ,d_\gamma -1}$ send it to
\[
y_r - y_m + 1-F(s,r) = 0. 
\]
Recalling that $F(s,r) = F(m,r)$, so it is exactly $H_{\varepsilon_r -\varepsilon_m,d_{\varepsilon_r-\varepsilon_m}}$.
Since $\kappa_{\varepsilon_r -\varepsilon_m} = F(r,m)+F(m,r)-1 = F(r,s) + F(s,r) -1 =1$, the resulting hyperplane belongs to $\mathcal{U}$.


The $y$-coordinate of a point belonging to a Scopes hyperplane satisfies: 
\[
\exists s\neq m\in \{ \pm 1,\ldots ,\pm e \}, \quad y_s -y_m \in \{ -F(s,m),-F(s,m)+1,\ldots , F(m,s) \}. 
\]
Assume that the class of $x\in \mathbb{Z}^e$ in $(\mathbb{Z}/q\mathbb{Z})^e$ avoids the Scopes hyperplanes and $-q/2 \le y_s \le q/2$ for any $s\in \{ 1,\ldots ,e\}$. 
If $x_s \ge x_{-s}$, then $F(-s,s)< y_s -y_{-s} < q-F(s,-s)$, which implies that $y_s \not\in \{ 0,q/2 \}$. 
Also, for $s,m$ such that $0< y_s \le y_m$, then $y_m - y_s > F(s,m)$. 
In particular, $y_s \neq y_m$ whenever $s\neq m\in \{\pm 1,\ldots ,\pm e\}$.
Similar to the proof of \cref{lem:gap-count}, a point $x\in (\mathbb{Z}/q\mathbb{Z})^e$ avoiding the Scopes hyperplane is uniquely determined by a signed permutation $\sigma = (\sigma_1,\ldots ,\sigma_e)$ and a sequence $0 < y_1 < \ldots < y_e < q/2$ such that
\begin{align*}
 2y_1&\geq1+F(\bar\sigma_1,\sigma_1),\\
 y_{i+1}-y_i&\geq1+F(\sigma_i,\sigma_{i+1}),\\
 q-2y_e&\geq1+F(\sigma_e,\bar\sigma_e),
\end{align*}
here we replace the strict triangle formula \eqref{eq:strict-triangle} by $1+F(s,m) < 2+F(s,r) + F(r,m)$.

Set $t=q-2e$ and $y_i'=y_i-i+\tfrac12$. The conditions become the following
\begin{equation}\label{eq:modified_avoidance}
\begin{aligned}
 2y_1'&\geq F(\bar\sigma_1,\sigma_1),\\
 y_{i+1}'-y_i'&\geq F(\sigma_i,\sigma_{i+1}),\\
 t-2y'_e&\geq F(\sigma_e,\bar\sigma_e).
\end{aligned}
\end{equation}
In particular $0\leq y_1'\leq\cdots\leq y_e'\leq t/2$.
The transformation is reversible once the signed order is retained.
The original label coset was $y_i\in-h_{\sigma_i}+\mathbb Z$; the new
one is $y_i'\in-h_{\sigma_i}+\tfrac12+\mathbb Z$. Undoing the sign and
gauge, namely setting
\[
 X_j=\operatorname{sgn}(\sigma_i)y_i'+h_j = x_j - \sigma^{-1}_j + \operatorname{sgn}(\sigma^{-1}_j)\frac{1}{2}\quad
          \text{when }|\sigma_i|=j,
\]
puts every coordinate in $\tfrac12+\mathbb Z$. 
Equivalently, $y'_i = X_{\sigma_i}-h_{\sigma_i}$.
The avoidance conditions become:
\begin{equation}\label{eq:shifted_criterion}
\sigma_m^{-1} < \sigma^{-1}_s \Rightarrow b_{\varepsilon_m -\varepsilon_s} (1) \le X_s - X_m \le  t -b_{\varepsilon_s-\varepsilon_m}(1).
\end{equation}
Let $\gamma =\varepsilon_s - \varepsilon_m$ with $s\neq m\in \{\pm1,\ldots ,\pm e\}$.
If $\kappa_\gamma =1$, then \eqref{eq:shifted_criterion} implies that $X_s -X_m \not\equiv -d_{-\gamma} = d_\gamma \pmod t$.
If $\kappa_\gamma =0$, then every residue modulus $t$ occurs once in the interval $[b_{-\gamma}(1),t-b_\gamma(1)]$.
If $\kappa_\gamma =-1$, every residue occurs once except for that of $b_{-\gamma}(1)$, which occurs twice. 
Since $2$ is invertible in $\mathbb{Z}/t\mathbb{Z}$, the map $z\mapsto z+t/2$ yields an $H$-equivariant isomorphism $((1/2+\mathbb{Z})/t\mathbb{Z})^e \to (\mathbb{Z}/t\mathbb{Z})^e$.
The class of the image of $X$ avoids $\mathcal{U}$.

Conversely, given a point $X'\in (\mathbb{Z}/t\mathbb{Z})^e$ that avoids $\mathcal{U}$, we have multiple choices of pairs $(X,\sigma)$ such that $X\equiv X'-\frac{t}{2}(1,\ldots ,1) (\operatorname{mod}\; t\mathbb{Z}^e)$ with $|X_s -h_s|\le t/2$ and \eqref{eq:shifted_criterion} holds. 
Set $(y'_1,\ldots ,y'_e)=\sigma^{-1}(X_1-h_1,\ldots ,X_e-h_e)$, then $0\le y'_1\le \ldots \le y'_e \le t/2$ and \eqref{eq:modified_avoidance} hold. 

1) In any admissible non-decreasing rearrangement of $(|X_1,-h_1|,\ldots ,|X_e-h_e|)$, the position of $|X_s-h_s|$ is smaller than $|X_m - h_m|$ if $|X_s-h_s| < |X_m - h_m|$ or $|X_s -h_s| = |X_m - h_m|$ and $F(s,m) >0$.
Recall that $F(s,m) + F(m,s) \ge 0$ and that the equality holds if and only if $\kappa_{\varepsilon_s-\varepsilon_m} =-1$. 
This determines the admissible choices of $(|\sigma_1|,\ldots ,|\sigma_e|)$.

2) The signs of $\sigma$ and $X$ are uniquely determined except for those $s$ with $X_s - h_s \in \{ 0,\pm t/2 \}$.
If $X_s - h_s =0$, then $\sigma^{-1}_s >0$ implies that $F(-s,s) =0$ and $\sigma^{-1}_s <0$ implies that $F(s,-s) =0$, multiple choices occur only when $\kappa_{2\varepsilon_s} =-1$.
Similarly, if $|X_s - h_s| =t/2$, multiple choices of $X_s$ and $\sigma^{-1}_t$ occur only when $\kappa_{2\varepsilon_s} =-1$.

Combining 1) and 2), there is a free $H_{X'}$-action on the set of these pairs $(X,\sigma)$ over $X'$.
In summary, the combination $x \mapsto (X,\sigma) \mapsto X'$ from $(\mathbb{Z}/q\mathbb{Z})^e\backslash \mathcal{A}_{a,1} \to (\mathbb{Z}/t\mathbb{Z})^e\backslash \mathcal{U}$ is surjective whose fiber at $X'$ is in bijection with $H_{X'}$.
So the number of points outside the Scopes arrangement equals to the left hand side of \eqref{eq:weighted-core-count} and its identification to the right hand side is obvious. 


All counts above use arbitrary odd integer moduli, not only fields.
Inclusion-exclusion and Smith normal form show that the original count
agrees with its characteristic polynomial for all sufficiently large
$q$ avoiding finitely many prime divisors. The same argument applies to
each fixed-point subsystem on the right at modulus $t$. There are
infinitely many simultaneously good pairs $q=t+2e$: at an odd prime at
most two residues for $t$ are excluded, while at $2$ the conditions
coincide since $2e$ is even. The Chinese remainder theorem supplies an
infinite progression. 
Agreement there proves the asserted polynomial
identity.
\end{proof}

\subsection{Evaluation of the zero-content core}

We record the elementary endpoint count needed below. 
For $1\leq u\leq v\leq n$, let
\[
\mathcal{J} (n,u,v) := \{ H_{\varepsilon_i - \varepsilon_j,0}: i<j \} \bigcup \{ H_{\varepsilon_i + \varepsilon_j,0}: j > \max (i,u)\} \bigcup \{ H_{\varepsilon_i,0} : i>v \}.
\]
Put $m=n-u$ and suppose $m\geq1$. With $(x)_k$ denoting a falling factorial, its complement over a large odd finite field $\mathbb{F}_t$ has size
\[
 (t-2m)_u\prod_{j=0}^{m-1}(t-1-2j)
 +(v-u)(t-2m+1)_u\prod_{j=0}^{m-2}(t-1-2j).
\]
Indeed the last $m$ coordinates are distinct up to sign; either none is
zero, or precisely one of the $v-u$ permitted labels is zero. Once these
are chosen, the first $u$ coordinates are chosen without replacement
from the remaining values. Factoring gives
\begin{equation}\label{eq:endpoint-core-factors}
 \chi(\mathcal J(n,u,v),t)=
 (t-(2n-v-1))^{1-\delta_{v,0}}\prod_{j=1}^{m}(t-(2j-1))
 \prod_{j=0}^{u-2}(t-2m-j)
 .
\end{equation}
Empty products are one and $\delta_{v,0}$ is the Kronecker delta. For $u=v=0$ the arrangement is the full $C_n$
arrangement, with factors $t-1,t-3,\ldots,t-(2n-1)$.

We define a signed component $\mathcal{J}_n$.   
As a hyperplane arrangement, $\mathcal{J}_n=\mathcal{J}(n,0,n)=\mathcal{J}(n,1,n)$ is the hyperplane arrangement of type $D_n$.
Let $H$ be the coordinate sign-change group, it's natural action preserves $\mathcal{J}_n$. 
Only the identity and single-coordinate flips can fix complement points. 
That is, $|H_X|\le 2$ for any complement point $X$ and the equality holds if and only if $X_i =0$ for some $i$.
Therefore
\begin{equation}\label{eq:signed-D-core}
 \sum_{X\notin \mathcal{J}(n,1,n)}|H_X|
 =\prod_{j=0}^{n-1}(t-1-2j)
       +2n\prod_{j=0}^{n-2}(t-1-2j)
 =(t+1)\prod_{j=0}^{n-2}(t-1-2j).
\end{equation}

\begin{theorem}\label{thm:zero-content-factor}
Let $\Lambda \in P^+$ in type $C^{(1)}_e$ and $\nu =\ell\Lambda_0 + a-\zeta_\Lambda (a)\delta \in P^+(\Lambda)$. 
Suppose that $\alpha=\Lambda-\nu =\sum c_i\alpha_i$ has a zero coefficient. 
Then there are integers $o_1,\ldots,o_e$, called offsets, such that for every $M\geq1$
\[
 \chi_{a,M}(q)
       =\prod_{i=1}^e(q-2eM-o_i).
\]
\end{theorem}
\begin{proof}

We apply \cref{thm:weighted-core} to compute the characteristic polynomial. 
To do this, we need to determine $\kappa_\gamma$'s.
By \cref{prop:zero-content-upward}, $N_\nu \subseteq \Phi_{J_1}\sqcup \ldots \sqcup \Phi_{J_k}$, where $J = \{i : c_i>0 \}$  and $J_i$'s are connected components of $J$.
Let $\gamma \in \mathring{\Phi}^+$. 
If $\gamma \in N_\nu$, then $d^+_\gamma = d^-_\gamma =0$ and $\kappa_\gamma =1$. 
If $\delta -\gamma\in N_\nu$, then $d^-_\gamma =-1$, $d^+_\gamma =1$ and thus $\kappa_\gamma =1$.
Otherwise, by \cref{thm:band-decomposition}, $b_\gamma (1)\in \{-1,0\}$ and is equal to $-1$ only if $\gamma = \varepsilon_i+\varepsilon_j$ and $a_i=a_j =\ell$. 
Similarly, $b_{-\gamma} (1)\in \{0,1\}$ and is equal to $0$ only if $(\gamma ,\nu) =0$. Hence, $\kappa_\gamma = -1$ if and only if $(\gamma ,\nu) \in \{0,2\ell\}$.

Combining with \cref{lem:zero-content-boundary}, we obtain a partition $I_0 =\{ 1,\ldots,e \} = I_1\sqcup I_2 \sqcup\ldots \sqcup I_c$: 
let $J'\subset I_0$ be the subset consisting of all $j$ such that $\{j,j- 1\}\cap J\neq \emptyset$, decompose $I_0\backslash J'$ according to the value of $a$.
Then $I_i$ is either a component of $J'$ or of $I_0 \backslash J'$.
Therefore, $(\mathbb{Z}/t\mathbb{Z})^e\backslash \mathcal{U} \cong (\mathbb{Z}/t\mathbb{Z})^{I_1}\backslash \mathcal{U}_{I_1} \times \cdots \times (\mathbb{Z}/t\mathbb{Z})^{I_c}\backslash \mathcal{U}_{I_c}$ and $H \cong H_1\times \cdots \times H_c$. 

(1) If $I_i = J_{s_i}$, then $e \in J_{s_i}$. 
Let $n=|I_i|$. 
As in \cref{prop:zero-content-upward}, we identify $I_i$ with $\{1,\ldots ,n \}$.
Then, $\mathcal{U}_{I_i} = \mathcal{J}(n,u,v)$ where $u =\max (\{j<n: b_j =1\}\cup \{ 0\})$ and $v = u$ unless $u=1$ in which case $v = \max (\{ 1\le j \le n : a_j = a_1  \})$.

(1-i) $u =0$, then $H_i$ is the trivial group, $\mathcal{U}_i$ is the full arrangement contributing offsets $\operatorname{Exp}(C_n)$.

(1-ii) $u=1$, then $H_i$ is trivial unless $a_1 =0$. 
When $H_i$ is trivial, the offsets are $1,3,\ldots ,2n-3,2n-v-1$.
When $H_i$ is non-trivial, it is generated by $H_{2\varepsilon_j,0}$ for all $1\le j \le n$, that is, the coordinate sign-change group. 
So the offsets are $-1,1,3,\ldots ,2n-3$.

(1-iii) $u>1$, then $H_i$ is trivial since $a_1 >0$. 
The offsets are $1,3,\ldots ,2m-1$, $2m,2m+1,\ldots , 2m+u-1$ with $m = n-u$.

(2) $I_i = J_{s_i} +1$, then $0 \in J_{s_i}$. We can deduce it from case (1) by applying the diagram involution $i\mapsto e-i$.

(3) $I_i = J_{s_i}\cup (J_{s_i}+1)$ with $0,e\notin J_{s_i}$.
The value of $a$ on $I_i$ lie in the open interval $(0,\ell)$ and all $\gamma = \varepsilon_{s}-\varepsilon_m$ supported in $I_i$ belong to $N_\nu$, so $H_i$ is trivial and $\mathcal{U}_{I_i}$ is the full $A_n$ arrangement, so the offsets are $0,1,\ldots ,n-1$. 

(4) $I_i\cap J' =\emptyset$ and $a_j \notin \{0,\ell\}$ for any $j\in I_i$. 
In this case, $H_i \cong \mathfrak{S}_n$ and $\mathcal{U}_{I_i} =\emptyset$. 
Then $\sum_{g\in H_i} ((\mathbb{Z}/t\mathbb{Z})^e)^g = t(t+1)\cdots (t+n-1)$.

(5) $I_i\cap J' =\emptyset$ and $a_j \in \{0,\ell\}$ for any $j\in I_i$. 
By symmetry, we may assume that $a$ takes value $0$ in $I_i$. 
In this case, $H_i$ is the group of signed permutation on $\{1,\ldots ,n\}$ and $\mathcal{U}_{I_i} = \emptyset$.
Then $\sum_{g\in H_i}(\mathbb{Z}/t\mathbb{Z})^g = \prod_{j=0}^{n-1}(t+1+2j)$.
The offsets are $-1,-3,\ldots ,1-2n$.
This completes the proof of $M =1$. The general result follows from \cref{thm:depth-shift}.
\end{proof}

\begin{theorem}\label{thm:all-level-product}
Let \(\Lambda=\sum_{j=0}^{\ell-1}\Lambda_{k_j}\) have positive
level in type \(C_e^{(1)}\), and let
\(\nu=\ell\Lambda_0+a-\zeta_\Lambda(a)\delta\) be a dominant
maximal supported weight. There is an explicit multiset
\(\mathcal E_{\Lambda,\nu}\) of \(e\) integers, independent of
\(M\), such that, for every integer \(M\ge1\),
\begin{equation}\label{eq:all-level-product}
 \chi_M(q)=\prod_{b\in\mathcal E_{\Lambda,\nu}}(q-2eM-b).
\end{equation}
If every coefficient of \(\Lambda-\nu\) is positive, put
\(r_A=L\), \(r_B=R-L\), \(r_C=e-R\), \(d=r_A+r_C\),
as in \cref{thm:positive-content-factor}. Then
\[
 \mathcal E_{\Lambda,\nu}
 =\{2d-u:u\in\operatorname{Exp}(C_{r_A})\sqcup
                       \operatorname{Exp}(C_{r_C})\}
      \sqcup\{2d,2d+1,\ldots,2d+r_B-1\}.
\]
If some coefficient is zero, the finite-component construction in
\cref{thm:zero-content-factor} specifies
\(\mathcal E_{\Lambda,\nu}\), including zero entries on unused
dimensions. In particular its computation requires only the content
coefficients and equal-coordinate blocks of the normalized dominant
vector.

Put \(n_j=\#\{i:a_i=j\}\). The number of Scopes classes in
\(W(\nu-M\delta)\) is the closed product
\begin{equation}\label{eq:all-level-scopes-product}
 \frac{\displaystyle\prod_{b\in\mathcal E_{\Lambda,\nu}}
                       (2eM+b+1)}
 {\displaystyle 2^{n_0+n_\ell}\prod_{j=0}^{\ell}n_j!}.
\end{equation}
\end{theorem}
\begin{proof}
The two alternatives for the content are exhaustive. Their
factorizations are \cref{thm:positive-content-factor} and
\cref{thm:zero-content-factor}. \cref{thm:all-chambers} and
the stabilizer calculation in \eqref{eq:all-chambers} give the product
count by evaluating at \(q=-1\).
\end{proof}

\begin{example}\label{ex:all-level-products}
For \(e=4\), \(\mathbf k=(0,3,4)\) and
\(a=(3,2,2,0)\), one has \(\zeta=1\) and
\(\Lambda-\nu=(1,1,1,1,1)\) in simple-root coordinates.
Thus \(L=1\), \(R=3\), and
\(\mathcal E_{\Lambda,\nu}=\{3,3,4,5\}\). Consequently
\[
 \chi_M(q)=(q-8M-3)^2(q-8M-4)(q-8M-5).
\]
The stabilizer order is eight, so there are
\(\frac18(8M+4)^2(8M+5)(8M+6)\) Scopes classes;
at depth one this is \(3276\).
For \(e=2\), \(\mathbf k=(0,1)\), \(a=(1,0)\),
the maximal content is zero and the stabilizer is
\(A_1\). Thus \(\mathcal E=\{-1,0\}\)
with stabilizer order two, and the counts are
\(2M(4M+1)\), namely \(10,36,78\) at depths one, two and three.
\end{example}

\section{The type-A comparison}\label{sec:typeA}

We prove the type-A product by ordinary circular orders. Thus its
characteristic-polynomial and class-count assertions require no
logarithmic freeness theorem.

\begin{theorem}[All-level type-$A$ factorization]
\label{thm:typeA-all-level-factorization}
Let $n\geq2$, let $\Lambda$ be a dominant integral weight
of $A_{n-1}^{(1)}$, and let $\nu\in P^+(\Lambda)$ be a maximal weight.
Write
\[
 \Lambda-\nu=\sum_{i\in\mathbb Z/n\mathbb Z}c_i\alpha_i,
 \qquad \nu_i=\langle\nu,\alpha_i^\vee\rangle,
 \qquad J=\{i:c_i>0\},\quad K=\{i:c_i=0,\ \nu_i=0\}.
\]
Define the multiset of $n-1$ integers
\begin{equation}\label{eq:typeA-offsets}
 \mathcal E=\operatorname{Exp}(W_J)
  \sqcup\{-d:d\in\operatorname{Exp}(W_K)\}
  \sqcup\{0\}^{\,n-1-|J|-|K|}.
\end{equation}
For every integer $M\geq1$,
\begin{equation}\label{eq:typeA-factorization}
 \boxed{\displaystyle
 \chi_{\Lambda,\nu,M}^{A_{n-1}}(q)
       =\prod_{d\in\mathcal E}(q-nM-d).}
\end{equation}
The number of Scopes classes in the Weyl orbit of $\nu-M\delta$ is
\begin{equation}\label{eq:typeA-class-count}
 \frac{1}{|W_\nu|}\prod_{d\in\mathcal E}(nM+d+1).
\end{equation}
If $I=\{i:\nu_i=0\}$ has components with $s_1,\ldots,s_t$ indexes,
then $|W_\nu|=\prod_{j=1}^t(s_j+1)!$.
\end{theorem}
\begin{proof}
Note that $\delta = \sum_{i=0}^{n-1}\alpha_i$. 
If $\min c_i >0$, then $\nu + \delta \in P^+$ and $\nu +\delta \le \Lambda$, thus $\nu +\delta\in P(\Lambda)$, a contradiction.
Therefore, $J$ is proper. 
Similar to \cref{prop:zero-content-upward}, $N_\nu = \Phi_J^+$.

If $c_i=0$ and $i\pm 1 \in J$, 
then $\nu_i =\Lambda_i-\sum_j a_{ij}c_j> 0$.
Hence, $i\notin K$. 
The subdiagram on $J\cup K$
is proper, and its components have mutually orthogonal gradient spans.
This includes $n=2$, using the doubled off-diagonal Cartan entries.
The band decomposition now gives
\begin{equation}\label{eq:typeA-signed-correction}
 \kappa=\mathbf1_{\mathcal A(\Phi_J)}
              -\mathbf1_{\mathcal A(\Phi_K)}.
\end{equation}

Use coordinate labels $1,\ldots,n$, modulo common translation, and put
\[
 d_{ij}=\min\{d\in\Z:\nu+\varepsilon_i-\varepsilon_j-d\delta \in P(\Lambda)\},\qquad f_{ij}=1-d_{ij}.
\]
By \cref{lem:absent-root-closure}, 
\begin{equation}\label{eq:typeA-weak-triangle}
 f_{ij}\leq f_{ik}+f_{kj},\qquad
 f_{ij}+f_{ji}=1+\kappa_{ij}\in\{0,1,2\}.
\end{equation}
Set $h_1=0$, $h_i=f_{1i}$ for $i\ne 1$, and write $F_{ij}=f_{ij}+h_i-h_j$. 
The triangle inequality and \eqref{eq:typeA-weak-triangle} imply $F_{ij}\geq0$, also when $i=1$ or $j=1$.
This is an integral change of coordinate origins.

Let $x = (x_1,\ldots ,x_n)\in \mathbb{Z}\mathring{\Phi} \subset \mathbb{Z}^n$ be a point whose class modulus $q$ lies in the complement of Scopes arrangement. 
For $i=1,\ldots ,n$, let $y'_i \in [0,q)$ be an integer such that $y'_i \equiv x_i-h_i+u\pmod q$ with $\sum h_i -nu \equiv 0\pmod q$. 
The avoidance conditions imply that there is a unique $\sigma \in\mathfrak{S}_n$ such that $0 \le y_1 \le \ldots \le y_n < q$ with $y_i = y'_{\sigma_i}$ and 
\[
\begin{aligned}
    & y_{i+1} - y_i \ge M + F_{\sigma_i,\sigma_{i+1}}\quad 1\le i < n \\
    & q - y_n +y_1\ge M +F_{\sigma_n,\sigma_1}
\end{aligned}
\]
Since $\sum y_i \equiv 0 \pmod q$, $y$ is uniquely determined once the differences $y_{i+1} - y_i$ are fixed.
Therefore, the number of complement point corresponding to $\sigma$ is
\begin{equation}\label{eq:typeA-circular-sum}
 \binom{q-nM-\sum_{i=1}^n F_{\sigma_i,\sigma_{i+1}}+n-1}{n-1},
 \qquad \sigma_{n+1}=\sigma_1.
\end{equation}
The proof of this identity follows a change of variables analogous to that used in the proof of \cref{lem:gap-count}.
\cref{lem:band-modular-count} applies for every sufficiently large
integer modulus in type A. 
Summing \eqref{eq:typeA-circular-sum} over
the symmetric group $\mathfrak{S}_n$ proves
$\chi_M(q)=\chi_1(q-n(M-1))$.

To evaluate the depth-one polynomial, 
we follow a strategy similar to \cref{thm:weighted-core}.
Define
\[
 \mathcal U=\{X_i-X_j=d_{ij}:\kappa_{ij}=1\},\qquad
 H=\langle s_{ij,1-d_{ji}}:\kappa_{ji}=-1\rangle.
\]
Here $s_{ij,c}$ is the reflection with respect to the hyperplane $X_i-X_j = d_{ij} -1$.
Given $(y,\sigma)$ as above, set $X_{\sigma_i} = y_i+h_{\sigma_i}-i +v$, where $v$ is a translation so that $\sum X_i \equiv 0\pmod q$. 
The avoidance conditions become
\[
\sigma^{-1}_i < \sigma^{-1}_j \Rightarrow f_{ij} \le X_j - X_i \le t-f_{ji},
\]
where $t =q-n$. 
Equivalently, the class of $X$ in $(\mathbb{Z}/t\mathbb{Z})^n$ belong to the complement of $\mathcal{U}$.
The fiber consisting of $(y,\sigma)$ over a class $X\notin \mathcal{U}$ is a $H_X$-torsor.
In summary,
\begin{equation}\label{eq:typeA-weighted-core}
 \chi_1(q)=\sum_{X\in (\mathbb{Z}/t\mathbb{Z}\mathring{\Phi})\backslash \mathcal{U}}|H_X| =\frac1t
    \sum_{X\in(\Z/t\Z)^n\setminus\mathcal U}|H_X|,
\end{equation}
for all sufficiently large $q$. 

The core hyperplanes of $\mathcal U$ are those of $\Phi_J$, and the
reflections of $H$ are those of $\Phi_K$ by \eqref{eq:typeA-signed-correction}.
By applying a suitable diagram automorphism, we may assume that $0 \notin J\cup K$ in which case the reflections correspond to $\mathcal{U}$ or belong to $H$ are all coordinate permutation.
The components of $J$ and $K$ determines a partition $[1,n] = I_1 \sqcup \ldots \sqcup I_c$ so that $n_i =|I_i| >1$ if and only if it is of the form $J_{c_i}\cup J_{c_i +1}$ (called positive block) for some component $J_{c_i}$ of $J$ or of the form $K_{c_i}\cup K_{c_i}+1$ (called negative block) for some component of $K$. 
Then, $(\mathbb{Z}/t\mathbb{Z})^n\backslash \mathcal{U} \cong \prod (\mathbb{Z}/t\mathbb{Z})^{n_i}\backslash \mathcal{U}_i$ is compatible with $H \cong \prod H_i$.
If $I_i$ is positive, then $H_i$ is trivial and $\mathcal{U}_i\neq \emptyset$ is a full arrangement of type $A_{n_i-1}$ which contributes a factor
\[
 t(t-1)\cdots(t-n_i+1)
\]
to the sum in right hand side \eqref{eq:typeA-weighted-core}. 
If $I_i$ is negative, then $H_i$ is non-trivial and $\mathcal{U}_i=\emptyset$, it contributes a factor
\[
 \sum_{\pi\in\mathfrak S_{n_i}}t^{\#\text{cycles}(\pi)}
       =t(t+1)\cdots(t+n_i-1).
\]
Other $I_i$'s contribute factors $t$ as $H_i=\{e\}$ and $(\mathbb{Z}/t\mathbb{Z})^{n_i}\backslash \mathcal{U}_i=\mathbb{Z}/t\mathbb{Z}$.
This shows \eqref{eq:typeA-offsets}.
The depth shift proves \eqref{eq:typeA-factorization}.
Finally, \cref{thm:all-chambers} and Zaslavsky's formula give
\eqref{eq:typeA-class-count}. The order of $W_\nu$ follows from its
finite type-A components.
\end{proof}

\begin{remark}[Classification and enumeration]\label{rem:typeA-comparison}
Chuang--Tan's Corollary~3.22 \cite{CT} classifies all Scopes classes
by their moving vectors, Scopes vectors and pyramid numbers. It thus
determines a total count by enumerating the realizable invariant
tuples. Their displayed product in Corollary~3.37 counts the RoCK
classes, rather than all classes. \cref{thm:typeA-all-level-factorization}
evaluates the total as a closed product in dominant maximal-weight
coordinates; it does not assert a new classification of type-A blocks.

For $A_2^{(1)}$ and $\nu=\Lambda$, the table compares the total with
Chuang--Tan's RoCK product $3!/\prod_i k_i!=6/|W_\Lambda|$.
The nonzero multiplicities $k_i$ in its three rows are $(3)$, $(2,1)$
and $(1,1,1)$, respectively. The RoCK count is independent of $M\geq1$.
\[
\begin{array}{c|c|ccc|c}
 \Lambda&\mathcal E& M=1&M=2&M=3&\text{RoCK}\\\hline
 \Lambda_0&\{-2,-1\}&1&5&14&1\\
 \Lambda_0+\Lambda_1&\{-1,0\}&6&21&45&3\\
 \Lambda_0+\Lambda_1+\Lambda_2&\{0,0\}&16&49&100&6
\end{array}
\]
\end{remark}

\subsection{A limit on extension to other untwisted types}

The highest-weight band theorem is uniform in type. The product
formula at an arbitrary maximal weight is not: already in
$G_2^{(1)}$ its characteristic polynomial can be irreducible over $\R$.

\begin{proposition}\label{prop:G2-counterexample}
In $G_2^{(1)}$, number the finite short root $\alpha_1$ and the finite
long root $\alpha_2$, so that
$\delta=\alpha_0+3\alpha_1+2\alpha_2$. Take
\[
 \Lambda=\Lambda_0+\Lambda_1,\qquad
 \nu=\Lambda-\alpha_0-\alpha_1-\alpha_2.
\]
Then $\nu$ is dominant and maximal. At depth one its Scopes
arrangement has characteristic polynomial
\[
 \chi_1(q)=q^2-18q+83
\]
and its orbit contains $17$ Scopes classes. In particular integral
factorization at every dominant maximal weight does not extend to
all untwisted affine types.
\end{proposition}
\begin{proof}
The Cartan matrix, with coroots indexing rows, is
\[
 \begin{pmatrix}2&0&-1\\0&2&-3\\-1&-1&2\end{pmatrix}.
\]
Thus $\nu$ has Dynkin labels $(0,2,0)$ and level two. 
$\nu\in P(\Lambda)$ by the dominant support criterion, whereas $\nu+\delta \ge \Lambda$.
The positive real roots $\le \Lambda-\nu$ are exactly
\[
 \alpha_0,\quad\alpha_1,\quad\alpha_2,\quad
 \alpha_0+\alpha_2,\quad\alpha_1+\alpha_2,\quad
 \alpha_0+\alpha_1+\alpha_2.
\]
All six belong to $N_\nu$. 
Indeed, this is a consequence of the following:
\[
\begin{aligned}
    & \Lambda - \alpha_0 = s_0 (\Lambda), \qquad \Lambda - \alpha_1 = s_1 (\Lambda), \\
    & \Lambda - \alpha_0-\alpha_1 = s_1s_0(\Lambda),\\
    & \Lambda - \alpha_0-\alpha_2 = s_2s_0(\Lambda),\\
    & \Lambda - \alpha_1-\alpha_2 = s_2s_1(\Lambda).
\end{aligned}
\]

Use finite coordinates $x=\alpha_1(z)$ and $y=\alpha_2(z)$.
The normalized point of $\nu$ is $(1/3,0)$.
\cref{thm:band-decomposition} gives precisely the eighteen lines
\begin{equation}\label{eq:G2-walls}
 \begin{aligned}
 x,\ y,\ x+y&\in\{-1,0,1\},\\
 2x+y,\ 3x+y,\ 3x+2y&\in\{0,1,2\}.
 \end{aligned}
\end{equation}
Solving pairs of these linear equations gives the following complete
intersection tally, where the multiplicity is the number of lines
through the point:
\[
 \begin{array}{c|rrrr}
  \text{multiplicity }m&2&3&4&6\\\hline
  \text{number of points}&33&13&3&3
 \end{array}
\]
For an independent check of completeness, these entries account for
$33+3\cdot13+6\cdot3+15\cdot3=135$ pairs, exactly
$\binom{18}{2}-6\binom32$, the number of nonparallel pairs.
The rank-two intersection formula now gives
\[
 \chi_1(q)=q^2-18q+(33+2\cdot13+3\cdot3+5\cdot3).
\]
Its discriminant is $-8$. The stabilizer generated by nodes $0,2$
has type $A_2$ and order six, so \cref{thm:all-chambers}
gives $\chi_1(-1)/6=102/6=17$ classes.
\end{proof}

\begin{remark}
Exact rational intersection counts in $G_2^{(1)}$ for all highest
weights and dominant maximal weights of levels $1$ through $5$, at
depths $1,2,3$, give $281$ configurations and $843$ polynomials.
Among them $42$ configurations fail integral factorization at every
tested depth. All satisfy
$\chi_M(q)=\chi_1(q-6(M-1))$ at those depths, and every highest-weight
case agrees with \cref{thm:affine-band-factorization}.
Thus these finite checks separate the depth-shift question from the
false uniform factorization assertion; they do not prove a general
depth-shift theorem outside types A and C.
\end{remark}

\section{Rouquier weights and RoCK categories}\label{sec:rock}

The following count uses Webster's general chamber mechanism.
Once every finite-root direction contributes a wall, there is one RoCK
chamber for each finite Weyl chamber; only the stabilizer quotient
remains. We include the direction argument and the type-C stabilizer
calculation. The type-A analogue is \cite[Corollary~3.37]{CT}.

\begin{theorem}\label{thm:rock-count}
Suppose \(M\ge1\), and put \(n_r=\#\{i:a_i=r\}\) for
\(0\le r\le\ell\). The number of RoCK Scopes classes in \(W\mu\) is
\begin{equation}\label{eq:rock-count}
 \frac{2^{\,e-n_0-n_\ell}e!}{\prod_{r=0}^{\ell}n_r!}.
\end{equation}
At level one this is \(\binom et\), where \(t=h(p)\).
At level one and depth zero there is just one RoCK Scopes class.
\end{theorem}
\begin{proof}
We first show that every finite-root direction occurs among the walls.
For \(\gamma\in\mathring{\Phi}^+\), both \(\gamma\) and
\(\delta-\gamma\) are positive real roots. 
For a dominant weight $\mu$, $(\gamma ,\mu)\ge 0$, $(\delta -\gamma ,\mu)\ge 0$ and their sum is positive.
Also, \(\mu+\delta \in P(\Lambda)\) 
because \(M\ge1\). 
Therefore, at least one of $\mu +\delta -\gamma$ and $\mu + \delta -(\delta -\gamma)$ is a weight of $V(\Lambda)$.
This gives a wall with direction \(\gamma\).

A finite arrangement with every finite-root direction has exactly one
RoCK chamber for each finite Weyl chamber: 
for each signed permutation $\sigma$, the chamber $\bigcap_{1\le i \le e} H^+_{b_{-\sigma (\alpha_i)}(M)\delta-\sigma(\alpha_i)}$ is a RoCK chamber, where $H^+_{\alpha} = \{x: (x,\alpha)>0\}$.
Indeed, there is no hyperplane intersecting this chamber by the weak triangle inequality.
This RoCK chamber is a translation of the finite Weyl chamber $\sigma (C_0)$.
Hence there are \(2^e e!\) such chambers before taking
the stabilizer quotient. The stabilizer of \(a\) has order
\[
 |W_\mu|=2^{n_0+n_\ell}\prod_{r=0}^{\ell}n_r!.
\]
Indeed, $W_\mu$ is the proper standard parabolic subgroup generated by $s_i$ with $(\mu,\alpha_i) =0$.
Taking orbits proves \eqref{eq:rock-count}. At level one,
\(a=u(t)\), giving the binomial coefficient. When the level and
depth are both as in the last assertion, \cref{thm:normal-forms}
gives one Scopes class; a finite arrangement always has a RoCK chamber,
so that class is RoCK.
\end{proof}

The positive-depth hypothesis in \eqref{eq:rock-count} is essential.
Depth zero does not mean \(N_\mu=\varnothing\) as we've seen in \cref{prop:positive-content-upward} and \cref{prop:zero-content-upward}.

\subsection{A paired-runner Rouquier condition}
We distinguish Webster's intrinsic notion of RoCK from the following
fixed-order abacus condition. 
At level one, let \(m\) be depth and \(t=h(p)\), and use \(L_{ab}\) from \eqref{eq:threshold}, also
when \(m=0\). Put \(\tau_{ab}=\max(0,L_{ab}-2)\).

\begin{definition}\label{def:rouquier}
A level-one weight of charge \(k\) is \emph{abacus Rouquier} 
at depth $m$
if
\begin{equation}\label{eq:rouquier}
 p_1\ge\max(1,m),\qquad
 p_{i+1}-p_i\ge\tau_{p_i,p_{i+1}}\quad(1\le i<e).
\end{equation}
At multicharge \(\mathbf k\), a weight is abacus Rouquier 
at depth $m$
if every component of every multipartition of that weight satisfies
\eqref{eq:rouquier}, at its own charge, paired vector and depth.
\end{definition}

This is a sufficient separation convention with a common runner order
in all components. In particular, it does not require each member of
a RoCK class to satisfy the inequalities. The bound one at depth zero
fixes the positive orientation. The higher-level quantifier is the
same kind of componentwise quantifier used for Rouquier weights in
the type-A comparison in \cite[Section~4.2.3]{Web23}.

\begin{remark}\label{rem:KL-thresholds}
The comparison with \cite[Section~4.1a]{KLivesey} can be made at
the level of inequalities. Their \(d\)-Rouquier \(\bar p\)-core
has runner counts \(b_1\ge d\) and
\(b_j-b_{j-1}\ge d-1\). Here, write
\(p_i=2q_i+\epsilon_i\). For \(m\ge1\), our inequalities become
\[
 q_1\ge\left\lceil\frac{m-\epsilon_1}{2}\right\rceil,
 \qquad
 q_{i+1}-q_i\ge
 \begin{cases}
 m-\epsilon_{i+1},&t<k,\\
 m-1+\delta_{\epsilon_i,1}\delta_{\epsilon_{i+1},0},&t=k,\\
 m-1+\epsilon_i,&t>k.
 \end{cases}
\]
Thus the internal gaps have the same baseline \(m-1\), with
corrections from folded parity, while the endpoint has a different
normalization. The quantities \(q_i\) are half paired differences,
not the runner counts of a \(\bar p\)-abacus. This is a comparison
of separation patterns, not an identification of their twisted
\(A_{2e}^{(2)}\) setting with our untwisted \(C_e^{(1)}\) setting.
\end{remark}

\begin{lemma}\label{lem:bad-roots}
At a level-one abacus Rouquier weight \(\nu\), 
\(\nu+\alpha\notin P(\Lambda)\) for the roots
\begin{equation}\label{eq:bad-roots}
 \begin{aligned}
 2\eps_i+n\delta&\qquad(1\le i\le e,\ n\ge0),\\
 \eps_i+\eps_j+n\delta&\qquad(i<j,\ n\ge0),\\
 \eps_j-\eps_i+n\delta&\qquad(i<j,\ n\ge1).
 \end{aligned}
\end{equation}
\end{lemma}
\begin{proof}
Write \(c_0=f_k(p)+m\). For a short root toggling coordinates
\(i,j\), put \(t'=t+2-2\epsilon_i-2\epsilon_j\) and
\(\eta=(|t'-k|-|t-k|)/2\). 
Directly from \(f_k\),
\begin{align*}
 f_k(p+2\eps_i)-f_k(p)&=p_i+1,\\
 f_k(p+\eps_i+\eps_j)-f_k(p)&=(p_i+p_j+1+\eta)/2,\\
 f_k(p-\eps_i+\eps_j)-f_k(p)&=(p_j-p_i+1+\eta)/2.
\end{align*}
Here \(-1\le\eta\le1\). The inequalities
\(\tau_{ab}+\tau_{bc}\ge\tau_{ac}\) show that the adjacent
gap bounds imply the same bound for every \(i<j\).
For \(m\ge1\), that bound, together with the parity of
\(p_j-p_i\), is exactly sufficient for
\[
 p_j-p_i\ge2m-1-\eta.
\]
Thus, the first increments $>m$ and the last $\ge m$. 
This implies that $f_k(p') > c_0 -n$. 
By \cref{thm:depth}, $\nu + \alpha \not\in P(\Lambda)$.
If \(m=0\), the vector is positive and weakly increasing.
The first two increments are at least one and the third is
nonnegative. This proves the same conclusion at depth zero.
\end{proof}

\subsection{An effective test at higher level}
We first rewrite the level-one inequalities without a depth parameter.
For \(x\in\mathcal P_k\), put \(t=h(x)\),
\(\epsilon_i=x_i\bmod2\), and
\[
 \kappa_{ab}^{k,t}=
 \begin{cases}2-a-b,&t<k,\\0,&t=k,\\a+b,&t>k.\end{cases}
\]
Define the \emph{Rouquier capacity} by \(\mathfrak r_k(x)=-1\)
if \(x_1<1\) or \(x\) is not weakly increasing, and otherwise by
\begin{equation}\label{eq:capacity}
 \mathfrak r_k(x)=\min\left\{x_1,
    \min_{1\le i<e}\left\lfloor
    \frac{x_{i+1}-x_i+2-\kappa_{\epsilon_i,\epsilon_{i+1}}^{k,t}}2
                     \right\rfloor\right\}.
\end{equation}
For $\mu\in P(\Lambda_k)$ with paired vector \(x\) and depth
\(m\), \cref{def:rouquier} is equivalent to
\(m\le\mathfrak r_k(x)\). Indeed, positivity and weak increase
give the depth-zero conditions; the remaining bounds are
\(m\le x_1\) and
\(2m\le x_{i+1}-x_i+2-\kappa_{\epsilon_i,\epsilon_{i+1}}^{k,t}\).

\begin{lemma}[Component bound]\label{lem:component-bound}
Assume that $\Lambda =\sum_{c=0}^{\ell-1}\Lambda_{k_c}\in P^+$, $\nu = \ell \Lambda_0 + P -c_0\delta \in P(\Lambda)$ 
with depth \(M\) and defect \(d\). 
Put \(K=\sum_c k_c\). 
For every multipartition of this content, with component vectors \(p^{(c)}\) and depths
\(m_c\), one has \(m_c\le M\) and
\begin{equation}\label{eq:variance}
 \sum_c\left\|p^{(c)}-\frac P\ell\right\|^2
 \le 4c_0+K-\frac{\norm P^2}{\ell}
 =\frac{2d+\ell K-\norm u^2}{\ell}=:V.
\end{equation}
If \(\ell\ge2\), each component satisfies the sharper bound
\begin{equation}\label{eq:component-ball}
 \left\|p^{(c)}-\frac P\ell\right\|^2
                  \le\frac{\ell-1}{\ell}V.
\end{equation}
Both bounds are uniform on a Weyl orbit, after recentering at \(P/\ell\).
\end{lemma}
\begin{proof}
The depth assertion follows from \eqref{eq:higher-depth}.
Since \(\sum_c f_{k_c}(p^{(c)})\le c_0\) and
\(4f_{k_c}(p^{(c)})\ge\norm{p^{(c)}}^2-k_c\), summing and
subtracting \(\norm P^2/\ell\) gives \eqref{eq:variance}.
The equality follows from \eqref{eq:higher-coordinates}.
Set \(v_c=p^{(c)}-P/\ell\), so \(\sum_c v_c=0\). Cauchy's
inequality gives
\(\sum_{b\ne c}\norm{v_b}^2\ge\norm{v_c}^2/(\ell-1)\),
which proves \eqref{eq:component-ball}. Finally defect is Weyl
invariant, so \(V\) is as well.
\end{proof}

\begin{theorem}[Exact finite Rouquier test]\label{thm:effective-rouquier}
Suppose \(\ell\ge2\), and let \(\ell \Lambda_0+P -c_0\delta \in P(\Lambda)\). Delete
the \(c\)-th charge to form \(\widehat{\mathbf k}_c\). For
\(x\in\mathcal P_{k_c}\), define
\begin{equation}\label{eq:component-max-depth}
 m_c^{\max}(x)=c_0-f_{k_c}(x)
                    -T_{\widehat{\mathbf k}_c}(P-x).
\end{equation}
Then the content is abacus Rouquier if and only if, for every \(c\)
and every \(x\in\mathcal P_{k_c}\) satisfying
\[
 \left\|x-\frac P\ell\right\|^2\le\frac{\ell-1}{\ell}V,
 \qquad m_c^{\max}(x)\ge0,
\]
one has
\begin{equation}\label{eq:effective-rouquier}
                 m_c^{\max}(x)\le\mathfrak r_{k_c}(x).
\end{equation}
Set \(R=\lceil\sqrt{(\ell-1)V/\ell}\rceil\). It suffices to
test at most \(\ell(2R+1)^e\) integer vectors, using
\(O(\ell(2R+1)^e e\log e)\) arithmetic operations.
\end{theorem}
\begin{proof}
Fix \(c,x\). 
By \cref{thm:higher}, $(\ell-1)\Lambda_0 + (P-x) -(c_0 -f_{k_c}(x))\delta \in P(\Lambda -\Lambda_{k_c})$ if and only if $c_0 - f_{k_c}(x) \ge T_{\widehat{\mathbf k}_c}(P-x)$.
Thus, \(x = p^{(c)}(\lambda^{(c)})\) for some $\blam$ if and only if \(m_c^{\max}(x)\ge0\). 
More precisely, its possible component depths range over the integers from zero to \(m_c^{\max}(x)\).
To realize any such depth, use a partition at vector \(x\) with
that depth, and use \cref{thm:higher} for the remaining components. 
In particular, the maximum is attained when the remaining components together have depth zero. 
This completes the proof of necessity.
The sufficiency follows from the equivalent description of Rouquier abacus.


The ball lies in the coordinate box
\(P_i/\ell-R\le x_i\le P_i/\ell+R\), which contains at most
\((2R+1)^e\) integer vectors. Filter by parity and the ball bound,
then evaluate \eqref{eq:component-max-depth} with
\cref{thm:closed-minimum}. Computing capacity costs \(O(e)\).
This proves the arithmetic-operation bound. It is not a polynomial
bound in the rank, and no enumeration of multipartitions is involved.
\end{proof}

\begin{example}\label{ex:capacity}
Let \(e=2\), \(\mathbf k=(2,2)\), \(P=(2,2)\). Then
\(T_{(2,2)}(P)=0\). For \(c_0=0\) or \(1\), the only possible
component vector is \((1,1)\), whose capacity is one. Hence both
contents are abacus Rouquier. For \(c_0\ge2\), taking
\(x=(1,1)\) gives \(m_c^{\max}(x)=c_0>1\), so the test fails.
For instance the multipartition \(((8),\varnothing)\) has content
\(2\delta\); its first component has vector \((1,1)\) and depth
two. The exact test therefore decides this whole content string:
precisely depths zero and one are abacus Rouquier in the chosen order.
\end{example}

\subsection{The Rouquier--RoCK comparison}

\begin{theorem}[Rouquier--RoCK comparison]\label{thm:rouquier-rock}
Fix \(\Lambda=\sum_c\Lambda_{k_c}\). For any integrable categorical
action of affine type \(C_e^{(1)}\) with support \(P(\Lambda)\):
\begin{enumerate}
\item every abacus Rouquier weight is RoCK;
\item every RoCK Scopes class contains an abacus Rouquier weight.
\end{enumerate}
Consequently a Scopes class is RoCK if and only if it contains an
abacus Rouquier weight. 
In particular the assertions apply to the minimal categorification by representations of cyclotomic KLR algebras
in \cref{sec:scopes}.
\end{theorem}
\begin{proof}
We first explain how to test RoCKness without making the weight
dominant. For \(\nu\in P(\Lambda)\), transport the arrangement
above to obtain
\begin{equation}\label{eq:transport-walls}
 \mathcal H_\nu=\{H_\alpha: \alpha \in \Phi_{\mathrm{re}}, \; \langle\nu,\alpha^\vee\rangle\ge0,\;
       \nu+\alpha\in P(\Lambda)\}.
\end{equation}
Here, \(\alpha\) ranges over real roots. 
This is precisely the transported set of walls: 
$\nu =w^{-1}(\mu)$ for $\mu\in P^+$, then $\mathcal{H}_\nu = w^{-1}(\mathcal{H}_\mu)$.
The chamber (with repsect to $\mathcal{H}_\nu$) containing \(\mathfrak a_\ell\) is RoCK exactly when \(\nu\) is RoCK.

Set \(D^- =\{0>x_1>\cdots>x_e\}\). 
The hyperplanes corresponding to the roots listed in \eqref{eq:bad-roots} are the only ones that can separate \(\mathfrak a_\ell\) from points far in \(D^-\). 
Therefore, to prove that Rouquier weight $\nu$ is RoCK, it suffices to prove, for every root in this list, both \(\langle\nu,\alpha^\vee\rangle\ge0\) and \(\nu+\alpha \notin P(\Lambda)\) as in this case, $\mathfrak{a}_\ell$ and $D^-$ lie in the same chamber. 

Suppose that \(\nu\) is abacus Rouquier. 
Let \(V_\otimes=\bigotimes_c V(\Lambda_{k_c})\). 
It is well-known that $P(V_{\otimes}) =P(\Lambda)$, see also 
\cref{prop:support}. 
Every decomposition \(\nu=\sum_c\nu_c\) in this tensor product is realized by a multipartitions. 
Hence, each \(\nu_c\) is Rouquier, and \cref{lem:bad-roots} gives
\((V(\Lambda_{k_c}))_{\nu_c+\alpha}=0\) for every root
\(\alpha\) in \eqref{eq:bad-roots}.
The primitive coproduct of a real-root raising operator therefore
makes \(E_\alpha\) zero on \((V_\otimes)_\nu\).

Each component paired vector is positive and weakly increasing, so
\(\langle\nu_c,\alpha^\vee\rangle\ge0\) for all the specified
roots, and hence \(\langle\nu,\alpha^\vee\rangle\ge0\).
Note that $V(\Lambda)_{\nu + \alpha}\neq 0$ implies that $E_\alpha|_{V(\Lambda)_\nu}$ is nonzero by the integrability of $V(\Lambda)$, a contradiction.
Thus, \(\nu+\alpha \notin P(\Lambda)\).
The preceding chamber argument proves (a).

For (b), fix a dominant \(\mu\in P(\Lambda)\) and a RoCK chamber \(R\).
Choose a signed permutation \(v\) such that \(R\) contains a
translate of \(vD^-\), and an integer vector \(z\) sufficiently
deep in that cone. With
\[
 w=t_{2\ell z}v,
 \qquad P=P(w^{-1}\mu)=v^{-1}a-2\ell v^{-1}z,
\] 
the alcove \(w\mathfrak a_\ell\) lies in \(R\). 
As \(z\) is taken deeper thus, $v^{-1}z\in D^-$ and \(P_1\) and $P_{i+1}-P_i$'s tend to $+\infty$. 
Note that $w\in W$, so \(w^{-1}\mu\) stays in the required orbit and class.

Define $V$ as in \eqref{eq:variance}, it depends only on $\Lambda$ and $\mu$.
Put $\sigma =\lceil\sqrt V\rceil$, and $z$ is deep enough so that
\begin{equation}\label{eq:stretch-bound}
 P_1\ge\ell\bigl(\max(1,M)+\sigma\bigr),\qquad
 P_{i+1}-P_i\ge2\ell(M+\sigma).
\end{equation}
For an arbitrary multipartition of this weight, we have, by \cref{lem:component-bound}, $|p_i^{(c)} - P_i/\ell| \le \sigma$ and $m_c\le M$.
In particular, $p^{(c)}_1 \ge \max (1,m_c)$ and $p^{(c)}_{i+1}-p_i^{(c)}\ge 2M\ge 2m_c$.
Hence, each component satisfies \eqref{eq:rouquier}. 
This completes the proof of (b).
\end{proof}

The theorem is the analogue of the two assertions in
\cite[Proposition~4.6]{Web23}, with \cref{def:rouquier}
specifying the type-C abacus condition. Its additional ingredients
are the paired-runner bounds and the uniform allocation estimate
\eqref{eq:variance}. It supplies abelian equivalences inside each
RoCK class, but does not identify these categories with a proposed
wreath product or another local algebra.

\begin{example}\label{ex:rouquier}
In \cref{ex:four-classes}, the two vectors \((1,2)\) and
\((2,3)\) satisfy \eqref{eq:rouquier}. At charge one and depth one,
they have partitions \(\lambda=(5,1,1)\) and \(\rho=(9,2,2)\) respectively:
\begin{center}
\begin{tabular}{cc}
\Tableau{12101,0,1}&\Tableau{121012101,01,10}\\[3pt]
\Abacus[runner labels={0,1,2,3},entries=beads]{4}{5,1,1,0,0}
&\Abacus[runner labels={0,1,2,3},entries=beads]{4}{9,2,2,0,0}\\
\(p=(1,2),\quad c=(2,4,1)\)&\(p=(2,3),\quad c=(4,7,2)\)
\end{tabular}
\end{center}
We display the runners by their $\ge -1$ parts.
Their normal forms are respectively
\((1,2)\) and \((-1,0)\). They therefore give two different RoCK
classes. \cref{thm:rock-count} gives exactly
\(\binom21=2\), so these are all the RoCK classes among the four
Scopes classes.
\end{example}

\begin{example}\label{ex:higher-rouquier}
Let \(e=2\), \(\Lambda=\Lambda_0+\Lambda_1\), and take the
orbit with \(a=(1,0)\) and source depth \(M=1\). Then
\(\zeta_\Lambda(a)=0\), \(d=4\), \(V=9/2\), and \(\sigma=3\).
The vector \(P=(9,28)\) belongs to this orbit and satisfies
\eqref{eq:stretch-bound}. Its least zeroth coefficient is
\[
 T_{(0,1)}(9,28)=\frac{9^2+28^2-1}{8}=108.
\]
Thus the content \((c_0,c_1,c_2)=(109,210,91)\), at source depth
one, is abacus Rouquier. \cref{thm:rock-count} gives four
RoCK classes in this orbit, since \(n_0=n_1=1\), \(n_2=0\).
\end{example}

\subsection{Labels and simple modules at depth zero}
The fibre count in \cref{prop:core-fibre} concerns
partition labels. The following refinement makes its difference from
a simple-module count explicit.

\begin{proposition}\label{prop:zero-simples}
At level one, let \(t=h(p)\) and \(D_k=0\), and define
\[
 (a,b)=
 \begin{cases}(t,(t-k)/2),&t\ge k,\\
 (e-t,(k-t)/2),&t\le k.
 \end{cases}
\]
The content has \(\binom ab\) partition labels, and its cyclotomic
KLR algebra has
\begin{equation}\label{eq:zero-simples}
 \binom ab-\binom a{b-1}
\end{equation}
simple modules, counted up to grading shift in the graded setting.
We use \(\binom a{-1}=0\).
\end{proposition}
\begin{proof}
The number of partitions is \eqref{eq:g}. 
Categorification identifies the number of simples with \(\dim V(\Lambda_k)_{\Lambda_0+p-f_k(p)\delta}\).
Weyl invariance lets us take \(p=u(t)\). Diagram reversal exchanges
\((k,t)\) with \((e-k,e-t)\), so assume \(t\le k\).
Then \(f_k(u(t))=0\), and the depth-zero weight is $\mu =\Lambda_0 + u(t) =\Lambda_t$.
As $\Lambda_k -\mu \in \mathbb{Z}\mathring{\Phi}$, $V(\Lambda_k)_{\mu}$ is the weight space of the $\mathfrak{sp}_{2e}$-submodule generated by a highest weight vector which is isomorphic to $V(\varpi_k)$.

Over characteristic zero, this is the primitive part of
\(\bigwedge^k\mathbb C^{2e}\); symplectic contraction onto
\(\bigwedge^{k-2}\mathbb C^{2e}\) is surjective for \(k\le e\).
At weight \(\varpi_t\), the first \(t\) positive basis vectors
are forced, and the other choices are \(b=(k-t)/2\) opposite
pairs from \(a=e-t\) indices. The two exterior-power weight
multiplicities are \(\binom ab\) and \(\binom a{b-1}\).
Their difference proves \eqref{eq:zero-simples}. The calculation is
in the characteristic-zero highest-weight module; categorification
gives the stated count over the field fixed in \cref{sec:scopes}.
\end{proof}

For \(|t-k|=2\), the defect is one. If \(t=k+2\), there are
\(k+2\) labels and \(k+1\) simples; if \(t=k-2\), there are
\(e-k+2\) labels and \(e-k+1\) simples. These agree with the
straight-line Brauer tree algebras in
\cite[Corollary~4.3 and Theorem~4.4]{CH}. Thus the term delta-core
in this paper means maximality along a \(\delta\)-string, and
includes positive-defect, nonsemisimple content summands.

\section*{Acknowledgements}
This project began in September 2025 at the University of Sydney, where both authors were supported by Australian Research Council Discovery Grant DP240101809.

We established the affine type \(C\) combinatorics and the level-one counting results before we began using AI assistance in September 2026. We subsequently used ChatGPT 6 Astra to assist in extending these results to higher levels. We have verified and simplified the arguments and revised the exposition for clarity. We take full responsibility for the mathematical correctness of the manuscript.

The first author acknowledges support from the National Natural Science Foundation of China (Grant No. 12301043).

We thank Andrew Mathas, Liron Speyer and Robert Muth for helpful conversations.

\bibliography{reference}
\bibliographystyle{alpha}

\end{document}